\documentclass[11pt]{article}
\usepackage{amsfonts}
\usepackage{graphics}
\usepackage{indentfirst}
\usepackage{color}
\usepackage{cite}
\usepackage{latexsym}
\usepackage[paper=a4paper, left=2.1cm, right=2.1cm, top=2.2cm, bottom=1.5cm, headheight=5.5pt, footskip=0.8cm, footnotesep=0.8cm, centering, includefoot]{geometry}
\usepackage{amsmath}
\allowdisplaybreaks
\usepackage{amssymb}
\usepackage[dvips]{epsfig}
\usepackage{amscd}

\newtheorem{theorem}{Theorem}[section]
\newtheorem{remark}{Remark}[section]
\newtheorem{definition}{Definition}[section]
\newtheorem{lemma}{Lemma}[section]

\DeclareMathOperator{\curl}{curl}
\DeclareMathOperator{\loc}{loc}

\DeclareMathOperator{\divv}{div}

\makeatletter
\@addtoreset{equation}{section}
\makeatother
\makeatletter
\@addtoreset{equation}{section}
\makeatother

\title{{Singularity formation to the Cauchy problem of the two-dimensional non-baratropic magnetohydrodynamic equations without heat conductivity}
\thanks{Supported by Chongqing
Research Program of Basic Research and Frontier Technology (No. cstc2018jcyjAX0049), the Postdoctoral Science Foundation of Chongqing (No. xm2017015), and China Postdoctoral Science Foundation (Nos. 2018T110936, 2017M610579).}
}

\author{Xin Zhong\thanks{School of Mathematics and Statistics, Southwest University, Chongqing 400715,
People's Republic of China ({\tt xzhong1014@amss.ac.cn}).
}
}
\date{ }

\begin{document}
\maketitle

\begin{abstract}
We study the singularity formation of strong solutions
to the two-dimensional (2D) Cauchy problem of the non-baratropic compressible magnetohydrodynamic equations without heat conductivity. It is proved that the strong solution exists globally if the density and the pressure are bounded from above. In particular, the criterion is independent of the magnetic field and is just the same as that of the compressible Navier-Stokes equations. Our method relies on weighted energy estimates and a Hardy-type inequality.
\end{abstract}


Math Subject Classification: 76W05; 35B65

\section{Introduction}
Let $\Omega\subset\mathbb{R}^2$ be a domain, the motion of a viscous, compressible, and heat conducting magnetohydrodynamic (MHD) flow in $\Omega$ can be described by the non-baratropic compressible MHD equations
\begin{align}\label{1.1}
\begin{cases}
\rho_{t}+\divv(\rho\mathbf{u})=0,\\
(\rho\mathbf{u})_{t}+\divv(\rho\mathbf{u}\otimes\mathbf{u})
-\mu\Delta\mathbf{u}
-(\lambda+\mu)\nabla\divv\mathbf{u}+\nabla P=
\mathbf{b}\cdot\nabla\mathbf{b}-\frac12\nabla|\mathbf{b}|^2,\\
c_{\nu}[(\rho\theta)_{t}+\divv(\rho\mathbf{u}\theta)]
+P\divv\mathbf{u}-\kappa\Delta\theta
=2\mu|\mathfrak{D}(\mathbf{u})|^2+\lambda(\divv\mathbf{u})^2
+\nu|\nabla\times\mathbf{b}|^2,\\
\mathbf{b}_{t}-\mathbf{b}\cdot\nabla\mathbf{u}
+\mathbf{u}\cdot\nabla\mathbf{b}+\mathbf{b}\divv\mathbf{u}=\nu\Delta\mathbf{b},\\
\divv \mathbf{b}=0.
\end{cases}
\end{align}
Here, $t\geq0$ is the time, $x\in\Omega$ is the spatial coordinate, and $\rho, \mathbf{u}, P=R\rho\theta\ (R>0), \theta, \mathbf{b}$ are the fluid density, velocity, pressure, absolute temperature, and the magnetic field respectively; $\mathfrak{D}(\mathbf{u})$ denotes the deformation tensor given by
\begin{equation*}
\mathfrak{D}(\mathbf{u})=\frac{1}{2}(\nabla\mathbf{u}+(\nabla\mathbf{u})^{tr}).
\end{equation*}
The constant viscosity coefficients $\mu$ and $\lambda$ satisfy the physical restrictions
\begin{equation}\label{1.2}
\mu>0,\ \mu+\lambda\geq0.
\end{equation}
Positive constants $c_\nu,\ \kappa$, and $\nu$ are respectively the heat capacity, the ratio of the heat conductivity coefficient over the heat capacity, and the magnetic diffusive coefficient.

There is huge literature on the studies about the theory of
well-posedness of solutions to the Cauchy problem and the initial boundary value problem (IBVP) for the compressible MHD system due to the physical importance, complexity, rich phenomena, and mathematical challenges, refer to \cite{HW2008,HW2010,LSX2016,LXZ2013,FY2009,HHPZ2017} and references therein.
However, many physical important and mathematical fundamental problems are still open due to the lack of smoothing mechanism and the strong nonlinearity. Kawashima \cite{K1983} first obtained the global existence and
uniqueness of classical solutions to the multi-dimensional compressible MHD equations when the initial data are close to a non-vacuum equilibrium in $H^3$-norm.
When the initial density allows vacuum, the local well-posedness of strong solutions to the initial boundary value problem of 3D non-isentropic MHD equations has been obtained by Fan-Yu \cite{FY2009}. For general large initial data, Hu-Wang \cite{HW2008,HW2010} proved the global existence of weak solutions with finite energy in Lions' framework for compressible Navier-Stokes equations \cite{L1998,F2004} provided the adiabatic exponent is suitably large. Recently, Li-Xu-Zhang \cite{LXZ2013} established the global existence and uniqueness of classical solutions to the Cauchy problem for the isentropic compressible MHD system in 3D with smooth initial data which are of small energy but possibly large oscillations and vacuum, which generalized the result for compressible Navier-Stokes equations obtained by Huang-Li-Xin \cite{HLX2012}. Very recently, Hong-Hou-Peng-Zhu \cite{HHPZ2017} improved the result in \cite{LXZ2013} to allow the initial energy large as long as the adiabatic exponent is close to 1 and $\nu$ is suitably large.
Furthermore, L{\"u}-Shi-Xu \cite{LSX2016} established the global existence and uniqueness of strong solutions to the 2D MHD equations provided that the smooth initial data are of small total energy. Nevertheless, it is an outstanding challenging open problem to investigate the global well-posedness for general large strong solutions with vacuum.

Therefore, it is important to study the mechanism of blow-up and structure of possible singularities of strong (or classical) solutions to the compressible MHD equations. The pioneering work can be traced to \cite{HX2005}, where He and Xin proved Serrin's criterion for strong solutions to the incompressible MHD system, that is,
\begin{equation}\label{1.3}
\lim_{T\rightarrow T^{*}}\|\mathbf{u}\|_{L^s(0,T;L^r)}=\infty,\ \text{for}\ \frac2s+\frac3r=1,\ 3<r\leq\infty,
\end{equation}
here $T^*$ is the finite blow up time. For the Cauchy problem of 2D compressible isentropic MHD system, Wang \cite{W2016} obtained the following criterion
\begin{equation}\label{1.4}
\lim_{T\rightarrow T^{*}}\|\rho\|_{L^\infty(0,T;L^\infty)}=\infty.
\end{equation}
This criterion asserts that the concentration of density must be responsible for the loss of regularity in finite time. For the IBVP of 2D full compressible MHD system,  Fan-Li-Nakamura \cite{FLN2015} proved that
\begin{equation}\label{1.5}
\lim_{T\rightarrow T^{*}}\left(\|\divv\mathbf{u}\|_{L^1(0,T;L^\infty)}
+\|\mathbf{b}\|_{L^\infty(0,T;L^\infty)}\right)
=\infty.
\end{equation}
Later on, Lu-Chen-Huang \cite{LCH2016} extended \eqref{1.5} with a refiner form
\begin{equation}\label{1.6}
\lim_{T\rightarrow T^{*}}\|\divv\mathbf{u}\|_{L^1(0,T;L^\infty)}
=\infty.
\end{equation}
The criterion \eqref{1.6} is the same as \cite{W2014} for 2D compressible full Navier-Stokes equations, which shows that the mechanism of blow-up is independent of the magnetic field.
Recently, for the Cauchy problem and the IBVP of 3D full compressible MHD system, Huang-Li \cite{HL2013} established the following Serrin type criterion
\begin{equation}\label{1.7}
\lim_{T\rightarrow T^{*}}\left(\|\rho\|_{L^\infty(0,T;L^\infty)}+\|\mathbf{u}\|_{L^s(0,T;L^r)}\right)
=\infty,\ \text{for}\ \frac2s+\frac3r\leq1,\ 3<r\leq\infty.
\end{equation}
There are also some interesting blow-up criteria for the compressible MHD system, see \cite{XZ2012,LDY2011}.


It should be noted that all the results mentioned above on the blow-up criteria of strong (or classical) solutions of viscous, compressible, and heat conducting MHD flows are for $\kappa>0$. Recently, for the 3D non-isentropic compressible Navier-Stokes equations with $\kappa=0$, Huang-Xin \cite{HX2016} showed that
\begin{equation}\label{1.06}
\lim_{T\rightarrow T^{*}}\left(\|\rho\|_{L^\infty(0,T;L^\infty)}
+\|\theta\|_{L^\infty(0,T;L^\infty)}\right)
=\infty
\end{equation}
under the assumption
\begin{equation}\label{1.07}
\mu>4\lambda.
\end{equation}
Later on, for the MHD flows, the author \cite{Z20172} obtained
 \begin{equation}\label{1.08}
\lim_{T\rightarrow T^{*}}\left(\|\mathfrak{D}(\mathbf{u})\|_{L^{1}(0,T;L^\infty)}
+\|P\|_{L^{\infty}(0,T;L^\infty)}\right)=\infty
\end{equation}
provided that
\begin{equation}\label{1.09}
3\mu>\lambda.
\end{equation}
It is worth mentioning that in a well-known paper \cite{X1998}, Xin considered non-isentropic compressible Navier-Stokes equations with $\kappa=0$ in multidimensional space, starting with a compactly supported initial density. He first proved that if the support of the density grows sublinearly in time and if the entropy is bounded from below then the solution cannot exist for all time. One key ingredient in the proof is a differential inequality on some integral functional (see \cite[Proposition 2.1]{X1998} for details). As an application, any smooth solution to the full compressible Navier-Stokes equations for polytropic fluids in the absence of heat conduction will blow up in finite time if the initial density is compactly supported.
Recently, based on the key observation that if initially a positive mass
is surrounded by a bounded vacuum region, then the time evolution remains uniformly bounded for all time, Xin-Yan \cite{XY2013} improved the blow-up results in \cite{X1998} by removing the assumptions that the initial density has compact support and the smooth solution has finite energy, but the initial data only has an isolated mass group.
Thus it seems very difficult to study globally smooth solutions of full compressible Navier-Stokes equations without heat conductivity in
multi-dimension, the same difficulty also arises in multi-dimensional MHD equations. These motivate us to study a blow-up criterion for the system \eqref{1.1} with zero heat conduction.
In fact, this is the main aim of this paper.

When $\kappa=0$, and without loss of generality, take $c_\nu=R=1$, the system \eqref{1.1} can be written as
\begin{align}\label{1.10}
\begin{cases}
\rho_{t}+\divv(\rho\mathbf{u})=0,\\
(\rho\mathbf{u})_{t}+\divv(\rho\mathbf{u}\otimes\mathbf{u})
-\mu\Delta\mathbf{u}
-(\lambda+\mu)\nabla\divv\mathbf{u}+\nabla P=
\mathbf{b}\cdot\nabla\mathbf{b}-\frac12\nabla|\mathbf{b}|^2,\\
P_{t}+\divv(P\mathbf{u})
+P\divv\mathbf{u}
=2\mu|\mathfrak{D}(\mathbf{u})|^2+\lambda(\divv\mathbf{u})^2
+\nu|\nabla\times\mathbf{b}|^2,\\
\mathbf{b}_{t}-\mathbf{b}\cdot\nabla\mathbf{u}
+\mathbf{u}\cdot\nabla\mathbf{b}+\mathbf{b}\divv\mathbf{u}=\nu\Delta\mathbf{b},\\
\divv \mathbf{b}=0.
\end{cases}
\end{align}
The present paper aims at giving a blow-up criterion of strong solutions to the Cauchy problem  of the system \eqref{1.10} with the initial condition
\begin{equation}\label{1.11}
(\rho,\rho\mathbf{u},P,\mathbf{b})(x,0)
=(\rho_0,\rho_0\mathbf{u}_0,P_0,\mathbf{b}_0)(x),\ \ x\in\mathbb{R}^2,
\end{equation}
and the far field behavior
\begin{equation}\label{1.12}
(\rho,\mathbf{u},P,\mathbf{b})(x,t)\rightarrow(0, \mathbf{0},0, \mathbf{0}),\ \text{as}\ |x|\rightarrow+\infty,\ t>0.
\end{equation}

Before stating our main result, we first explain the notations and conventions used throughout this paper. For $r>0$, set
\begin{equation*}
B_r \triangleq\left.\left\{x\in\mathbb{R}^2\right|\,|x|<r \right\},
\quad \int \cdot dx\triangleq\int_{\mathbb{R}^2}\cdot dx.
\end{equation*}
For $1\leq p\leq\infty$ and integer $k\geq0$, the standard Sobolev spaces are denoted by:
\begin{equation*}
L^p=L^p(\mathbb{R}^2),\ W^{k,p}=W^{k,p}(\mathbb{R}^2), \ H^{k}=H^{k,2}(\mathbb{R}^2),
\ D^{k,p}=\{u\in L_{\loc}^1|\nabla^k u\in L^p\}.
\end{equation*}
Now we define precisely what we mean by strong solutions to the problem \eqref{1.10}--\eqref{1.12}.
\begin{definition}[Strong solutions]\label{def1}
$(\rho,\mathbf{u},P,\mathbf{b})$ is called a strong solution to \eqref{1.10}--\eqref{1.12} in $\mathbb{R}^2\times(0,T)$, if for some $q_0>2$
and $a>1$,
\begin{equation*}
\begin{split}
\begin{cases}
\rho\geq0,\ \rho\bar{x}^a\in C([0,T];L^1\cap H^1\cap W^{1,q_0}),\ \rho_t\in C([0,T];L^{q_0}),\\
(\mathbf{u},\mathbf{b})\in C([0,T];D^{1,2}\cap D^{2,2})\cap L^{2}(0,T;D^{2,q_0}),\
\mathbf{b}\in C([0,T];H^2), \\
(\mathbf{u}_t,\mathbf{b}_t)\in L^{2}(0,T;D^{1,2}),\
(\sqrt{\rho}\mathbf{u}_{t},\mathbf{b}_t)\in L^{\infty}(0,T;L^{2}), \\
P\geq0,\ P\in C([0,T];L^1\cap H^1\cap W^{1,q_0}),\ P_t\in C([0,T];L^{q_0}), \\
\end{cases}
\end{split}
\end{equation*}
and $(\rho,\mathbf{u},P,\mathbf{b})$ satisfies both \eqref{1.10} almost everywhere in $\mathbb{R}^2\times(0,T)$ and \eqref{1.11} almost everywhere in $\mathbb{R}^2$. Here
\begin{equation}\label{2.01}
\bar{x}\triangleq(e+|x|^2)^{\frac12}\log^{1+\eta_0}(e+|x|^2)
\end{equation}
and $\eta_0$ is a positive number.
\end{definition}

Without loss of generality, we assume that the initial density $\rho_0$ satisfies
\begin{equation}\label{1.8}
\int_{\mathbb{R}^2} \rho_0dx=1,
\end{equation}
which implies that there exists a positive constant $N_0$ such that
\begin{equation}\label{1.9}
\int_{B_{N_0}}  \rho_0  dx\ge \frac12\int\rho_0dx=\frac12.
\end{equation}

Our main result reads as follows:
\begin{theorem}\label{thm1.1}
In addition to \eqref{1.8} and \eqref{1.9}, assume that the initial data $(\rho_0\geq0, \mathbf{u}_0,P_0\geq0,\mathbf{b}_0)$ satisfies for any given numbers $a>1$ and $q>2$,
\begin{align}\label{A}
\begin{cases}
\rho_{0}\bar{x}^{a}\in L^{1}\cap H^{1}\cap W^{1,q},\ \sqrt{\rho_0}\mathbf{u}_0\in L^2 ,\ \nabla\mathbf{u}_{0}\in H^1,\\
P_0\in L^1\cap H^1\cap W^{1,q},\ \mathbf{b}_0\bar{x}^{\frac{a}{2}}\in H^{1},\ \nabla^2\mathbf{b}_0\in L^{2},\ \divv\mathbf{b}_0=0,
\end{cases}
\end{align}
and the compatibility conditions
\begin{equation}\label{A2}
-\mu\Delta\mathbf{u}_0-(\lambda+\mu)\nabla\divv\mathbf{u}_0
+\nabla P_0-\mathbf{b}_0\cdot\nabla\mathbf{b}_0+\frac12\nabla|\mathbf{b}_0|^2
=\sqrt{\rho_0}\mathbf{g}
\end{equation}
for some $\mathbf{g}\in L^2(\Omega)$.
Let $(\rho,\mathbf{u},P,\mathbf{b})$ be a strong solution to the problem \eqref{1.10}--\eqref{1.12}. If $T^{*}<\infty$ is the maximal time of existence for that solution,
then we have
\begin{align}\label{B}
\lim_{T\rightarrow T^{*}}\left(\|\rho\|_{L^\infty(0,T;L^\infty)}
+\|P\|_{L^{\infty}(0,T;L^\infty)}\right)=\infty.
\end{align}
\end{theorem}

\begin{remark}\label{re1.1}
The local existence of a strong solution with initial data as in Theorem \ref{thm1.1} was established in \cite{LH2016}. Hence, the maximal time $T^{*}$ is well-defined.
\end{remark}

\begin{remark}\label{re1.2}
It is worth noting that the blow-up criterion \eqref{B} is independent of the magnetic field, and it is just the same as that of compressible Navier-Stokes equations \cite{Z2017}. Thus we generalize \cite[Theorem 1.1]{Z2017} to the compressible MHD flows.
\end{remark}

\begin{remark}\label{re1.3}
Compared with \cite{Z20172}, where the author investigated a blow-up criterion for the 3D Cauchy problem of non-isentropic magnetohydrodynamic equations with zero heat conduction, there is no need to impose additional restrictions on the viscosity coefficients $\mu$ and $\lambda$ except the physical restrictions \eqref{1.2}.
\end{remark}

We now make some comments on the analysis of this paper. We mainly make use of continuation argument to prove Theorem \ref{thm1.1}. That is, suppose that \eqref{B} were false, i.e.,
\begin{equation*}
\lim_{T\rightarrow T^*}\left(\|\rho\|_{L^\infty(0,T;L^\infty)}
+\|P\|_{L^{\infty}(0,T;L^\infty)}\right)\leq M_0<\infty.
\end{equation*}
We want to show that
\begin{equation*}
\sup_{0\leq t\leq T^*}\left(\|(\rho,P)\|_{H^1\cap W^{1,q}}
+\|\rho\bar{x}^{a}\|_{L^{1}\cap H^{1}\cap W^{1,q}}
+\|\nabla\mathbf{u}\|_{H^1}+\|\mathbf{b}\|_{H^2}
+\|\mathbf{b}\bar{x}^{a}\|_{H^{1}}
\right) \leq C<+\infty.
\end{equation*}

It should be pointed out that the crucial techniques of proofs in \cite{FLN2015,LCH2016} cannot be adapted directly to the situation treated here, since their arguments depend crucially on the boundedness of the domains and $\kappa>0$. Moreover, technically, it is hard to modify the three-dimensional analysis of \cite{Z20172} to the two-dimensional case
with initial density containing vacuum since the
analysis of \cite{Z20172} depends crucially on the a priori $L^6$-bound on the velocity, while in two dimensions it seems difficult to bound the $L^p(\mathbb{R}^2)$-norm of $\mathbf{u}$ just in terms of $\|\sqrt{\rho}\mathbf{u}\|_{L^{2}(\mathbb{R}^2)}$ and $\|\nabla\mathbf{u}\|_{L^{2}(\mathbb{R}^2)}$ for any $p\geq1$.

To overcome these difficulties mentioned above, some new ideas are needed.
Inspired by \cite{LX2013,Z2017}, we first observe that if the initial density decays not too slow at infinity, i.e., $\rho_0\bar{x}^a\in L^1(\mathbb{R}^2)$ for some positive constant $a>1$ (see \eqref{A}), then for any $\eta\in(0,1]$, we can show that (see \eqref{06.2})
\begin{equation}\label{1.20}
\mathbf{u}\bar{x}^{-\eta}\in L^{p_0}(\mathbb{R}^2),\  \text{for some}\  p_0>1.
\end{equation}
To finish the
higher order estimates, our new observation is to obtain the $L_t^\infty L_x^2$-norm of $\bar{x}^{\frac{a}{2}}\mathbf{b}$ and $\bar{x}^{\frac{a}{2}}\nabla\mathbf{b}$ (see
Lemma \ref{lem40}).
Then, motivated by the technique of Hoff \cite{H1995}, in order to get the $L_t^\infty L_x^{2}$-norm of $\sqrt{\rho}\dot{\mathbf{u}}$, we first show the desired a priori estimates of the $L_t^\infty L_x^2$-norm of $\nabla\mathbf{u}$ and $\nabla\mathbf{b}$, which is the second key observation in this paper (see Lemma \ref{lem35}).
The a priori estimates on the $L_t^\infty L_x^{q}$-norm of $(\nabla\rho,\nabla P)$ can be obtained (see Lemma \ref{lem37}) by solving a logarithm Gronwall inequality based on a logarithm estimate for the Lam{\'e} system.
Finally, with the help of \eqref{1.20}, we can get the spatial weighted estimate of the density (see Lemma \ref{lem38}).

The rest of this paper is organized as follows. In Section \ref{sec2}, we collect some elementary facts and inequalities that will be used later. Section \ref{sec3} is devoted to the proof of Theorem \ref{thm1.1}.

\section{Preliminaries}\label{sec2}

In this section, we will recall some known facts and elementary inequalities that will be used frequently later.


We begin with the following Gagliardo-Nirenberg inequality (see \cite{N1959}).
\begin{lemma}[Gagliardo-Nirenberg]\label{lem22}
For $p\in[2,\infty), r\in(2,\infty)$, and $s\in(1,\infty)$, there exists some generic constant $C>0$ which may depend on $p,$ $r$, and $s$ such that for $f\in H^{1}(\mathbb{R}^2)$ and $g\in L^{s}(\mathbb{R}^2)\cap D^{1,r}(\mathbb{R}^2)$, we have
\begin{eqnarray*}
& & \|f\|_{L^p(\mathbb{R}^2)}^{p}\leq C\|f\|_{L^2(\mathbb{R}^2)}^{2}\|\nabla f\|_{L^2(\mathbb{R}^2)}^{p-2}, \\
& & \|g\|_{C(\overline{\mathbb{R}^2})}\leq C\|g\|_{L^s(\mathbb{R}^2)}^{s(r-2)/(2r+s(r-2))}\|\nabla g\|_{L^r(\mathbb{R}^2)}^{2r/(2r+s(r-2))}.
\end{eqnarray*}
\end{lemma}

The following weighted $L^m$ bounds for elements of the Hilbert space $\tilde{D}^{1,2}(\mathbb{R}^2)\triangleq\{ v \in H_{\loc}^{1}(\mathbb{R}^2)|\nabla v \in L^{2}(\mathbb{R}^2)\}$ can be found in \cite[Theorem B.1]{L1996}.
\begin{lemma}\label{1leo}
For $m\in [2,\infty)$ and $\theta\in (1+m/2,\infty),$ there exists a positive constant $C$ such that for all $v\in  \tilde{D}^{1,2}(\mathbb{R}^2),$
\begin{equation}\label{3h}
\left(\int_{\mathbb{R}^2} \frac{|v|^m}{e+|x|^2}\left(\log \left(e+|x|^2\right)\right)^{-\theta}dx  \right)^{1/m}
\leq C\|v\|_{L^2(B_1)}+C\|\nabla v\|_{L^2(\mathbb{R}^2) }.
\end{equation}
\end{lemma}

The combination of Lemma \ref{1leo} and the Poincar\'e inequality yields
the following useful results on weighted bounds, whose proof can be found in \cite[Lemma 2.4]{LX2013}.

\begin{lemma}\label{lem26}
Let $\bar x$ be as in \eqref{2.01}. Assume that $\rho \in L^1(\mathbb{R}^2)\cap L^\infty(\mathbb{R}^2)$ is a non-negative function such that
\begin{equation*} \label{2.i2}
\|\rho\|_{L^1(B_{N_1})} \ge M_1, \quad \|\rho\|_{L^1(\mathbb{R}^2)\cap L^\infty(\mathbb{R}^2)}\le M_2,
\end{equation*}
for positive constants $M_1, M_2$, and $ N_1\ge 1$. Then for $\varepsilon> 0$ and $\eta>0,$ there is a positive constant $C$ depending only on $\varepsilon,\eta, M_1,M_2$, and $ N_1$, such that every $v\in \tilde{D}^{1,2}(\mathbb{R}^2)$ satisfies
\begin{equation}\label{22}
\|v\bar x^{-\eta}\|_{L^{(2+\varepsilon)/\tilde{\eta}}(\mathbb{R}^2)}
\le C \|{\sqrt\rho}v\|_{L^2(\mathbb{R}^2)}+C \|\nabla v\|_{L^2(\mathbb{R}^2)},
\end{equation}
with $\tilde{\eta}=\min\{1,\eta\}$.
\end{lemma}


Next, for $\nabla^{\bot}\triangleq(-\partial_2,\partial_1)$, denoting the material derivative of $f$ by $\dot{f}\triangleq f_t+\mathbf{u}\cdot\nabla f$,
then we have the following $L^p$-estimate (see \cite[Lemma 2.5]{LSX2016}) for the elliptic system derived from the momentum equations \eqref{1.10}$_2$:
\begin{equation}\label{3.014}
\Delta F=\divv(\rho\dot{\mathbf{u}}-\divv(\mathbf{b}\otimes\mathbf{b})),\
\mu\Delta\omega=\nabla^{\bot}\cdot(\rho\dot{\mathbf{u}}-\divv(\mathbf{b}\otimes\mathbf{b})),
\end{equation}
where $F$ is the effective viscous flux, $\omega$ is vorticity given by
\begin{equation}\label{3.14}
F=(\lambda+2\mu)\divv\mathbf{u}-P-\frac12|\mathbf{b}|^2,\
\omega=\partial_1u_2-\partial_2u_1.
\end{equation}
\begin{lemma}\label{lem23}
Let $(\rho, \mathbf{u}, P, \mathbf{b})$ be a smooth solution of \eqref{1.10}. Then for $p\geq2$ there exists a positive constant $C$ depending only on $p, \mu$ and $\lambda$ such that
\begin{align}
&\|\nabla F\|_{L^p}+\|\nabla\omega\|_{L^p}
\leq C\left(\|\rho\dot{\mathbf{u}}\|_{L^p}
+\||\mathbf{b}||\nabla\mathbf{b}|\|_{L^p}\right), \label{2.3} \\
&\|F\|_{L^p}+\|\omega\|_{L^p}
\leq C\left(\|\rho\dot{\mathbf{u}}\|_{L^2}
+\||\mathbf{b}||\nabla\mathbf{b}|\|_{L^2}\right)^{1-\frac2p}
\left(\|\nabla\mathbf{u}\|_{L^2}+\|P\|_{L^2}+\|\mathbf{b}\|_{L^4}^2\right)^{\frac2p}, \label{2.4} \\
&\|\nabla\mathbf{u}\|_{L^p}
\leq C\left(\|\rho\dot{\mathbf{u}}\|_{L^2}
+\||\mathbf{b}||\nabla\mathbf{b}|\|_{L^2}\right)^{1-\frac2p}
\left(\|\nabla\mathbf{u}\|_{L^2}+\|P\|_{L^2}+\|\mathbf{b}\|_{L^4}^2\right)^{\frac2p}
+C\|P\|_{L^p}+\|\mathbf{b}\|_{L^{2p}}^2. \label{2.5}
\end{align}
\end{lemma}

Finally, the following Beale-Kato-Majda type inequality (see \cite[Lemma 2.3]{HLX20112}) will be used to estimate $\|(\nabla\rho,\nabla P)\|_{L^q}\ (q>2)$.
\begin{lemma}\label{lem25}
For $q\in(2,\infty)$, there is a constant $C(q)>0$ such that for all $\nabla\mathbf{v}\in L^2\cap D^{1,q}$, it holds that
\begin{equation}\label{2.2}
\|\nabla\mathbf{v}\|_{L^\infty}\leq C\left(\|\divv\mathbf{v}\|_{L^\infty}
+\|\curl\mathbf{v}\|_{L^\infty}\right)\log(e+\|\nabla^2\mathbf{v}\|_{L^q})
+C\|\nabla\mathbf{v}\|_{L^2}+C.
\end{equation}
\end{lemma}

\section{Proof of Theorem \ref{thm1.1}}\label{sec3}

Let $(\rho,\mathbf{u},P,\mathbf{b})$ be a strong solution described in Theorem \ref{thm1.1}. Suppose that \eqref{B} were false, that is, there exists a constant $M_0>0$ such that
\begin{equation}\label{3.1}
\lim_{T\rightarrow T^*}\left(\|\rho\|_{L^\infty(0,T;L^\infty)}
+\|P\|_{L^{\infty}(0,T;L^\infty)}\right)\leq M_0<\infty.
\end{equation}


First, we have the following standard estimate.
\begin{lemma}\label{lem32}
Under the condition \eqref{3.1}, it holds that for any $T\in[0,T^*)$,
\begin{equation}\label{3.3}
\sup_{0\leq t\leq T}\left(\|\sqrt{\rho}\mathbf{u}\|_{L^2}^2+\|\mathbf{b}\|_{L^2}^2+\|P\|_{L^1\cap L^\infty}
\right)+\int_{0}^T\left(\|\nabla\mathbf{u}\|_{L^2}^2
+\|\nabla\mathbf{b}\|_{L^2}^2\right)dt
\leq C,
\end{equation}
where and in what follows, $C,C_1,C_2$ stand for generic positive constants depending only on $M_0,\lambda,\mu,\nu,T^{*}$, and the initial data.
\end{lemma}
{\it Proof.}
1. It follows from \eqref{1.10}$_3$ that
\begin{equation}\label{7.01}
P_t+\mathbf{u}\cdot\nabla P+2P\divv\mathbf{u}
=F\triangleq2\mu|\mathfrak{D}(\mathbf{u})|^2
+\lambda(\divv\mathbf{u})^2+\nu|\nabla\times\mathbf{b}|^2\geq0.
\end{equation}
Define particle path before blowup time
\begin{align*}
\begin{cases}
\frac{d}{dt}\mathbf{X}(x,t)
=\mathbf{u}(\mathbf{X}(x,t),t),\\
\mathbf{X}(x,0)=x.
\end{cases}
\end{align*}
Thus, along particle path, we obtain from \eqref{7.01} that
\begin{align*}
\frac{d}{dt}P(\mathbf{X}(x,t),t)
=-2P\divv\mathbf{u}+F,
\end{align*}
which implies
\begin{equation}\label{3.03}
P(\mathbf{X}(x,t),t)=\exp\left(-2\int_{0}^t\divv\mathbf{u}ds\right)
\left[P_0+\int_{0}^t\exp\left(2\int_{0}^s\divv\mathbf{u}d\tau\right)Fds\right]\geq0.
\end{equation}

2. Multiplying \eqref{1.10}$_2$ and \eqref{1.10}$_3$ by $\mathbf{u}$ and $\mathbf{b}$, respectively, then adding the two resulting equations together, and integrating over $\mathbb{R}^2$, we obtain after integrating by parts that
\begin{align}\label{3.4}
\frac{1}{2}\frac{d}{dt}\int\left(\rho|\mathbf{u}|^2
+|\mathbf{b}|^2\right)dx
+\int\left[\mu|\nabla\mathbf{u}|^2+(\lambda+\mu)(\divv\mathbf{u})^2
+\nu|\nabla\mathbf{b}|^2\right]dx
=\int P\divv\mathbf{u}dx.
\end{align}
Integrating \eqref{1.10}$_3$ with respect to $x$ and then adding the resulting equality to \eqref{3.4} give rise to
\begin{align}\label{3.004}
\frac{d}{dt}\int\left(\frac12\rho|\mathbf{u}|^2+\frac12|\mathbf{b}|^2+P\right)dx=0,
\end{align}
which combined with \eqref{3.03}, \eqref{A}, and \eqref{3.1} leads to
\begin{align}\label{3.04}
\sup_{0\leq t\leq T}\left(\|\sqrt{\rho}\mathbf{u}\|_{L^2}^2+\|\mathbf{b}\|_{L^2}^2
+\|P\|_{L^1\cap L^\infty}
\right)\leq C.
\end{align}
This together with \eqref{3.4} and Cauchy-Schwarz inequality yields
\begin{align}\label{3.003}
\frac{d}{dt}\left(\|\sqrt{\rho}\mathbf{u}\|_{L^2}^2+\|\mathbf{b}\|_{L^2}^2\right)
+\mu\|\nabla\mathbf{u}\|_{L^2}^2+\nu\|\nabla\mathbf{b}\|_{L^2}^2
\leq C.
\end{align}
So the desired \eqref{3.3} follows from \eqref{3.04} and \eqref{3.003} integrated with respect to $t$. This completes the proof of Lemma \ref{lem32}.
\hfill $\Box$

Inspired by \cite{HX2005},
we have the following higher integrability of the magnetic field $\mathbf{b}$.
\begin{lemma}\label{lem33}
Under the condition \eqref{3.1}, it holds that
for any $p\geq2$ and $T\in[0,T^*)$,
\begin{equation}\label{3.5}
\sup_{0\leq t\leq T}\|\mathbf{b}\|_{L^p}+\int_0^T\int|\mathbf{b}|^{p-2}
|\nabla\mathbf{b}|^2dxdt\leq C.
\end{equation}
\end{lemma}
{\it Proof.}
Multiplying \eqref{1.10}$_4$ by $p|\mathbf{b}|^{p-2}\mathbf{b}$ and integrating the resulting equation over $\mathbb{R}^2$, we derive
\begin{align}\label{3.6}
\frac{d}{dt}\int|\mathbf{b}|^pdx
+\nu p(p-1)\int|\mathbf{b}|^{p-2}|\nabla\mathbf{b}|^2dx
 =p\int\left(\mathbf{b}\cdot\nabla\mathbf{u}-\mathbf{u}\cdot\nabla\mathbf{b}
-\mathbf{b}\divv\mathbf{u}\right)\cdot|\mathbf{b}|^{p-2}\mathbf{b}dx.
\end{align}
By the divergence theorem and \eqref{1.10}$_5$, we get
\begin{equation*}
-p\int(\mathbf{u}\cdot\nabla)\mathbf{b}\cdot|\mathbf{b}|^{p-2}\mathbf{b}dx
=\int\divv\mathbf{u}|\mathbf{b}|^pdx,
\end{equation*}
which together with \eqref{3.6} and Gagliardo-Nirenberg inequality yields
\begin{align*}
\frac{d}{dt}\int|\mathbf{b}|^pdx
+\nu p(p-1)\int|\mathbf{b}|^{p-2}|\nabla\mathbf{b}|^2dx
& \leq C\int|\nabla\mathbf{u}||\mathbf{b}|^pdx \notag \\
& \leq C\|\nabla\mathbf{u}\|_{L^2}\||\mathbf{b}|^{\frac{p}{2}}\|_{L^4}^2 \notag \\
& \leq C\|\nabla\mathbf{u}\|_{L^2}\||\mathbf{b}|^{\frac{p}{2}}\|_{L^2}
\|\nabla|\mathbf{b}|^{\frac{p}{2}}\|_{L^2}
 \notag \\
& \leq \varepsilon\|\nabla|\mathbf{b}|^{\frac{p}{2}}\|_{L^2}^2
+C\|\nabla\mathbf{u}\|_{L^2}^2\|\mathbf{b}\|_{L^p}^p,
\end{align*}
which along with Gronwall's inequality and \eqref{3.3} yields the desired \eqref{3.5} after Choosing $\varepsilon$ suitably small. This finishes the proof of Lemma \ref{lem33}. \hfill $\Box$

The following lemma gives the estimates on the spatial gradients of both the velocity and the magnetic field, which are crucial for deriving the higher order estimates of the solution.
\begin{lemma}\label{lem34}
Under the condition \eqref{3.1}, it holds that for any $T\in[0,T^*)$,
\begin{equation}\label{5.1}
\sup_{0\leq t\leq T}\left(\|\nabla\mathbf{u}\|_{L^2}^{2}
+\|\nabla\mathbf{b}\|_{L^2}^{2}\right)
+\int_{0}^{T}\left(\|\sqrt{\rho}\dot{\mathbf{u}}\|_{L^2}^{2}
+\|\nabla^2\mathbf{b}\|_{L^2}^2\right)dt \leq C.
\end{equation}
\end{lemma}
{\it Proof.}
1. Multiplying \eqref{1.10}$_2$ by $\dot{\mathbf{u}}$ and integrating the resulting equation over $\mathbb{R}^2$ give rise to
\begin{align}\label{5.2}
\int\rho|\dot{\mathbf{u}}|^2dx
& = -\int\dot{\mathbf{u}}\cdot\nabla Pdx+\mu\int\dot{\mathbf{u}}\cdot\Delta\mathbf{u}dx
+(\lambda+\mu)\int\dot{\mathbf{u}}\cdot\nabla\divv\mathbf{u}dx \notag \\
& \quad +\int\dot{\mathbf{u}}\cdot\mathbf{b}\cdot\nabla\mathbf{b}dx
-\frac12\int\dot{\mathbf{u}}\cdot\nabla|\mathbf{b}|^2dx\triangleq\sum_{i=1}^5I_i.
\end{align}
By \eqref{1.10}$_3$ and integrating by parts, we derive from \eqref{3.1} and Garliardo-Nirenberg inequality that
\begin{align}\label{5.3}
I_1
& = \int \left[(\divv\mathbf{u})_tP
-(\mathbf{u}\cdot\nabla\mathbf{u})\cdot\nabla P\right]dx \nonumber\\
&=\frac{d}{dt}\int P\divv\mathbf{u}dx+
\int\left[P(\divv\mathbf{u})^2-2\mu\divv\mathbf{u}|\mathfrak{D}(\mathbf{u})|^2
-\lambda(\divv\mathbf{u})^3-\nu\divv\mathbf{u}|\nabla\times\mathbf{b}|^2
\right]dx \nonumber\\
& \quad
+\int P\partial_ju_i\partial_iu_jdx \nonumber\\
&\leq \frac{d}{dt}\int P\divv\mathbf{u}dx+
C\|\nabla\mathbf{u}\|_{L^2}^2+C\|\nabla\mathbf{u}\|_{L^3}^3
+C\|\nabla\mathbf{b}\|_{L^3}^3\nonumber\\
&\leq \frac{d}{dt}\int P\divv\mathbf{u}dx+
C\|\nabla\mathbf{u}\|_{L^2}^2+C\|\nabla\mathbf{u}\|_{L^3}^3
+C\|\nabla\mathbf{b}\|_{L^2}^4+\frac{\nu}{4}\|\nabla^2\mathbf{b}\|_{L^2}^2.
\end{align}
It follows from integration by parts that
\begin{align}\label{5.4}
I_2
& =\mu\int(\mathbf{u}_{t} +\mathbf{u}\cdot\nabla\mathbf{u})\cdot\Delta\mathbf{u}dx \notag \\
& = -\frac{\mu}{2}\frac{d}{dt}\|\nabla\mathbf{u}\|_{L^2}^{2}
-\mu\int\partial_{i}u_{j}\partial_{i}(u_{k}\partial_{k}u_{j})dx \notag \\
 &  \leq-\frac{\mu}{2}\frac{d}{dt}\|\nabla\mathbf{u}\|_{L^2}^{2}
 +C\|\nabla\mathbf{u} \|_{L^3}^{3}.
\end{align}
Similarly to $I_2$, one gets
\begin{align}\label{5.5}
I_3
& = -\frac{\lambda+\mu}{2}\frac{d}{dt}\|\divv\mathbf{u}\|_{L^2}^{2}
-(\lambda+\mu)\int\divv\mathbf{u}\divv(\mathbf{u}\cdot\nabla\mathbf{u})dx \notag \\
 &  \leq-\frac{\lambda+\mu}{2}\frac{d}{dt}\|\divv\mathbf{u}\|_{L^2}^{2}
 +C\|\nabla\mathbf{u} \|_{L^3}^{3}.
\end{align}
By virtue of \eqref{1.10}$_4$ and \eqref{1.10}$_5$, one deduces from integration by parts and Gagliardo-Nirenberg inequality that
\begin{align}\label{5.6}
I_4
& = \int \mathbf{b}\cdot\nabla \mathbf{b}\cdot\mathbf{u}_tdx+\int\mathbf{b}\cdot\nabla\mathbf{b}\cdot
(\mathbf{u}\cdot \nabla\mathbf{u})dx
\notag \\
& = -\frac{d}{dt}\int\mathbf{b}\cdot\nabla\mathbf{u}\cdot\mathbf{b}dx
+\int\left(\mathbf{b}_t\cdot\nabla\mathbf{u}\cdot\mathbf{b}
+\mathbf{b}\cdot\nabla\mathbf{u}\cdot\mathbf{b}_t\right)dx
- \int\mathbf{b}\cdot\nabla(\mathbf{u}\cdot \nabla\mathbf{u})\cdot\mathbf{b}dx
\notag \\
& = -\frac{d}{dt}\int\mathbf{b}\cdot\nabla\mathbf{u}\cdot\mathbf{b}dx
+\int\left(\mathbf{b}_t\cdot\nabla\mathbf{u}\cdot\mathbf{b}
+\mathbf{b}\cdot\nabla\mathbf{u}\cdot\mathbf{b}_t\right)dx
- \int\left(\mathbf{b}\cdot\nabla u_i\partial_i\mathbf{u}\cdot \mathbf{b}+b_ku_i\partial_{ik}\mathbf{u}\cdot\mathbf{b}\right)dx
\notag \\
& = -\frac{d}{dt}\int\mathbf{b}\cdot\nabla\mathbf{u}\cdot\mathbf{b}dx
+\int\left(\mathbf{b}_t\cdot\nabla\mathbf{u}\cdot\mathbf{b}
+\mathbf{b}\cdot\nabla\mathbf{u}\cdot\mathbf{b}_t\right)dx \notag \\
& \quad - \int\left[\mathbf{b}\cdot\nabla u_k\partial_k\mathbf{u}\cdot \mathbf{b}-\mathbf{u}\cdot\nabla b_k\partial_{k}\mathbf{u}\cdot\mathbf{b}
-\mathbf{b}\cdot\nabla\mathbf{u}\cdot\mathbf{b}\divv\mathbf{u}-
(\mathbf{b}\cdot\nabla\mathbf{u})(\mathbf{u}\cdot\nabla\mathbf{b})\right]dx
\notag \\
& = -\frac{d}{dt}\int\mathbf{b}\cdot\nabla\mathbf{u}\cdot\mathbf{b}dx
+\int\left(\partial_tb_k-\mathbf{b}\nabla u_k+\mathbf{u}\cdot\nabla b_k\right)
\partial_k\mathbf{u}\cdot\mathbf{b}dx \notag \\
& \quad + \int\mathbf{b}\cdot\nabla\mathbf{u}\cdot(\mathbf{b}_t+\mathbf{b}\divv\mathbf{u}
+\mathbf{u}\cdot\nabla\mathbf{b})dx
\notag \\
& = -\frac{d}{dt}\int\mathbf{b}\cdot\nabla\mathbf{u}\cdot\mathbf{b}dx
+\int(\nu\Delta\mathbf{b}-\mathbf{b}\divv\mathbf{u})\cdot\nabla\mathbf{u}\cdot\mathbf{b}dx
+\int\mathbf{b}\cdot\nabla\mathbf{u}
\cdot(\nu\Delta\mathbf{b}+\mathbf{b}\cdot\nabla\mathbf{u})dx
\notag \\
& \leq -\frac{d}{dt}\int\mathbf{b}\cdot\nabla\mathbf{u}\cdot\mathbf{b}dx
+C\int\left(|\Delta\mathbf{b}||\nabla\mathbf{u}||\mathbf{b}|
+|\nabla\mathbf{u}|^2|\mathbf{b}|^2\right)dx
\notag \\
& \leq -\frac{d}{dt}\int\mathbf{b}\cdot\nabla\mathbf{u}\cdot\mathbf{b}dx
+C\|\nabla\mathbf{u}\|_{L^3}^3
+C\|\mathbf{b}\|_{L^6}^6+\frac{\nu}{8}\|\Delta\mathbf{b}\|_{L^2}^2
\notag \\
& \leq -\frac{d}{dt}\int\mathbf{b}\cdot\nabla\mathbf{u}\cdot\mathbf{b}dx
+C\|\nabla\mathbf{u}\|_{L^3}^3
+C\|\nabla\mathbf{b}\|_{L^2}^4+\frac{\nu}{4}\|\Delta\mathbf{b}\|_{L^2}^2.
\end{align}
Applying \eqref{1.10}$_4$, \eqref{1.10}$_5$, and Gagliardo-Nirenberg inequality, we have
\begin{align}\label{5.7}
I_5
& =\frac12\int|\mathbf{b}|^2\divv\mathbf{u}_tdx
+\frac12\int|\mathbf{b}|^2\divv(\mathbf{u}\cdot\nabla\mathbf{u})dx
\notag \\
& = \frac12\frac{d}{dt}\int|\mathbf{b}|^2\divv\mathbf{u}dx
-\frac12\int|\mathbf{b}|^2(\divv\mathbf{u})^2dx
+\frac12\int|\mathbf{b}|^2\partial_iu_j\partial_ju_idx \notag \\
& \quad -\int(\mathbf{b}\cdot\nabla\mathbf{u}+\nu\Delta\mathbf{b}-\mathbf{b}\divv\mathbf{u})
\cdot\mathbf{b}\divv\mathbf{u}dx \notag \\
& \leq \frac12\frac{d}{dt}\int|\mathbf{b}|^2\divv\mathbf{u}dx
+C\int|\mathbf{b}|^2|\nabla\mathbf{u}|^2dx
+\frac{\nu}{8}\|\Delta\mathbf{b}\|_{L^2}^2 \notag \\
& \leq \frac12\frac{d}{dt}\int|\mathbf{b}|^2\divv\mathbf{u}dx
+C\|\nabla\mathbf{u}\|_{L^3}^3+C\|\mathbf{b}\|_{L^6}^6
+\frac{\nu}{8}\|\Delta\mathbf{b}\|_{L^2}^2 \notag \\
& \leq \frac12\frac{d}{dt}\int|\mathbf{b}|^2\divv\mathbf{u}dx
+C\|\nabla\mathbf{u}\|_{L^3}^3+C\|\nabla\mathbf{b}\|_{L^2}^4
+\frac{\nu}{4}\|\Delta\mathbf{b}\|_{L^2}^2.
\end{align}
Putting \eqref{5.3}--\eqref{5.7} into \eqref{5.2}, we obtain from \eqref{2.5} and \eqref{3.1} that
\begin{align}\label{5.8}
\Psi'(t)+\|\sqrt{\rho}\dot{\mathbf{u}}\|_{L^2}^2
\leq C\|\nabla\mathbf{u}\|_{L^2}^2+C\|\nabla\mathbf{u}\|_{L^3}^3
+C\|\nabla\mathbf{b}\|_{L^2}^4+\frac{3\nu}{4}\|\Delta\mathbf{b}\|_{L^2}^2,
\end{align}
where
\begin{equation*}
\Psi(t)\triangleq\frac{\mu}{2}\|\nabla\mathbf{u}\|_{L^2}^{2}
+\frac{\lambda+\mu}{2}\|\divv\mathbf{u}\|_{L^2}^{2}
-\int P\divv\mathbf{u}dx+\int\mathbf{b}\cdot\nabla\mathbf{u}\cdot\mathbf{b}dx
-\frac12\int|\mathbf{b}|^2\divv\mathbf{u}dx
\end{equation*}
satisfies
\begin{equation}\label{5.9}
\frac{\mu}{2}\|\nabla\mathbf{u}\|_{L^2}^2-C
\leq \Psi(t)\leq\mu\|\nabla\mathbf{u}\|_{L^2}^2+C
\end{equation}
due to \eqref{3.3} and \eqref{3.5}.

2. Multiplying \eqref{1.10}$_4$ by $\Delta\mathbf{b}$ and integrating by parts, one deduces from H{\"o}lder's inequlity, Gagliardo-Nirenberg inequality, and \eqref{3.3} that
\begin{align}\label{5.10}
\frac{d}{dt}\int|\nabla\mathbf{b}|^2dx+2\nu\int|\Delta\mathbf{b}|^2dx
& \leq C\int|\nabla\mathbf{u}||\nabla\mathbf{b}|^2dx
+C\int|\nabla\mathbf{u}||\mathbf{b}||\Delta\mathbf{b}|dx \notag \\
& \leq C\|\nabla\mathbf{u}\|_{L^3}\|\nabla\mathbf{b}\|_{L^2}^{\frac43}
\|\nabla\mathbf{b}\|_{H^1}^{\frac23}
+C\|\nabla\mathbf{u}\|_{L^3}\|\mathbf{b}\|_{L^6}
\|\Delta\mathbf{b}\|_{L^2} \notag \\
& \leq  C\|\nabla\mathbf{u}\|_{L^3}^3+C\|\nabla\mathbf{b}\|_{L^2}^4
+\frac{\nu}{4}\|\Delta\mathbf{b}\|_{L^2}^2.
\end{align}
Adding \eqref{5.10} to \eqref{5.8}, we then derive from \eqref{2.5} and \eqref{3.1} that \begin{align}\label{5.11}
& B'(t)+\|\sqrt{\rho}\dot{\mathbf{u}}\|_{L^2}^2+\nu\|\Delta\mathbf{b}\|_{L^2}^2\notag \\
& \leq C\|\nabla\mathbf{u}\|_{L^2}^2+C\|\nabla\mathbf{u}\|_{L^3}^3
+C\|\nabla\mathbf{b}\|_{L^2}^4 \notag \\
& \leq C\|\nabla\mathbf{u}\|_{L^2}^2
+C\left(\|\rho\dot{\mathbf{u}}\|_{L^2}
+\||\mathbf{b}||\nabla\mathbf{b}|\|_{L^2}\right)\left(\|\nabla\mathbf{u}\|_{L^2}+1\right)^2
+C+C\|\nabla\mathbf{b}\|_{L^2}^4
 \notag \\
& \leq \frac12\|\sqrt{\rho}\dot{\mathbf{u}}\|_{L^2}^2
+C(1+\|\nabla\mathbf{u}\|_{L^2}^2+\|\nabla\mathbf{b}\|_{L^2}^2
+\||\mathbf{b}||\nabla\mathbf{b}|\|_{L^2}^2)
\left(1+\|\nabla\mathbf{u}\|_{L^2}^2+\|\nabla\mathbf{b}\|_{L^2}^2\right),
\end{align}
where
\begin{equation}\label{5.12}
\frac{\mu}{2}\|\nabla\mathbf{u}\|_{L^2}^2+\|\nabla\mathbf{b}\|_{L^2}^2-C
\leq B(t)\triangleq\Psi(t)+\|\nabla\mathbf{b}\|_{L^2}^2
\leq\mu\|\nabla\mathbf{u}\|_{L^2}^2+\|\nabla\mathbf{b}\|_{L^2}^2-C
\end{equation}
due to \eqref{5.9}.
Thus the desired \eqref{5.1} follows from \eqref{5.11}, \eqref{5.12}, \eqref{3.3}, \eqref{3.5}, and Gronwall's inequality. This completes the proof of Lemma \ref{lem34}.
\hfill $\Box$

The following spatial weighted estimate on the density showed in \cite[Lemma 3.5]{Z2017} plays a crucial role in deriving the higher order derivatives of the solutions $(\rho,\mathbf{u},P,\mathbf{b})$, we sketch it here for completeness.
\begin{lemma}\label{lem36}
Under the condition \eqref{3.1}, it holds that for any $T\in[0,T^*)$,
\begin{equation}\label{06.1}
\sup_{0\leq t\leq T}\|\rho\bar{x}^{a}\|_{L^{1}}\leq C(T).
\end{equation}
\end{lemma}
{\it Proof.}
1. For $N>1,$ let $\varphi_N\in C^\infty_0(\mathbb{R}^2)$  satisfy
\begin{equation} \label{vp1}
0\le \varphi_N \le 1, \quad  \varphi_N(x)
=\begin{cases} 1,~~~~ |x|\le N/2,\\
0,~~~~ |x|\ge N,\end{cases}
\quad |\nabla \varphi_N|\le C N^{-1}.
\end{equation}
It follows from \eqref{1.10}$_1$ that
\begin{align}\label{oo0}
\frac{d}{dt}\int \rho \varphi_{N} dx &=\int \rho\mathbf{u}
\cdot\nabla \varphi_{N} dx \notag \\
&\ge - C N^{-1}\left(\int\rho dx\right)^{1/2}
\left(\int\rho |u|^2dx\right)^{1/2}\ge - \tilde{C} N^{-1},
\end{align}
where in the last inequality one has used \eqref{3.3} and
\begin{equation*}
\int\rho dx=\int \rho_0dx.
\end{equation*}
Integrating \eqref{oo0} and choosing $N=N_1\triangleq2N_0+4\tilde CT$, we obtain after using \eqref{1.9} that
\begin{align}\label{p1}
\inf\limits_{0\le t\le T}\int_{B_{N_1}} \rho dx&\ge \inf\limits_{0\le t\le T}\int \rho \varphi_{N_1} dx \notag \\
&\ge \int \rho_0 \varphi_{N_1} dx-\tilde{C}N_1^{-1}T \notag \\
&\ge \int_{B_{N_0}} \rho_0 dx-\frac{\tilde{C}T}{2N_0+4\tilde{C} T} \notag \\
&\ge 1/4.
\end{align}
Hence, it follows from \eqref{p1}, \eqref{3.1}, \eqref{22}, \eqref{3.3}, and \eqref{5.1} that for any $\eta\in(0,1]$ and any $s>2$,
\begin{equation}\label{06.2}
\|\mathbf{u}\bar{x}^{-\eta}\|_{L^{s/\eta}}\leq C\left(\|\sqrt{\rho}\mathbf{u}\|_{L^2}+\|\nabla\mathbf{u}\|_{L^2}\right)
\le C.
\end{equation}

2. Multiplying \eqref{1.10}$_1$ by $\bar{x}^{a}$ and integrating the resulting equality by parts over $\mathbb{R}^2$ yield that
\begin{equation*}
\begin{split}
\frac{d}{dt}\int\rho\bar{x}^{a}dx & \leq C\int\rho|\mathbf{u} |\bar{x}^{a-1}\log^{2}(e+|x|^2)dx\\
 &  \leq C\|\rho\bar{x}^{a-1+\frac{8}{8+a}}\|_{L^{\frac{8+a}{7+a}}}\|\mathbf{u} \bar{x}^{-\frac{4}{8+a}}\|_{L^{8+a}} \\
 &  \leq C\int\rho\bar{x}^{a}dx+C,
\end{split}
\end{equation*}
which along with Gronwall's inequality gives \eqref{06.1} and finishes the proof of Lemma \ref{lem36}.    \hfill $\Box$


\begin{lemma}\label{lem40}
Under the condition \eqref{3.1}, and let $a>1$ be as in Theorem \ref{thm1.1}, then it holds that for any $T\in[0,T^*)$,
\begin{align}
\sup_{0\leq t\leq T}\|\mathbf{b}\bar{x}^{\frac{a}{2}}\|_{L^2}^2
+\int_{0}^T\|\nabla\mathbf{b}\bar{x}^{\frac{a}{2}}\|_{L^2}^2dt\leq C,
\label{10.1} \\
\sup_{0\leq t\leq T}\|\nabla\mathbf{b}\bar{x}^{\frac{a}{2}}\|_{L^2}^2
+\int_{0}^T\|\nabla^2\mathbf{b}\bar{x}^{\frac{a}{2}}\|_{L^2}^2dt\leq C. \label{10.2}
\end{align}
\end{lemma}
{\it Proof.}
1. Multiplying \eqref{1.10}$_4$ by $\mathbf{b}\bar{x}^a$ and integrating by parts give rise to
\begin{align}\label{10.3}
\frac12\frac{d}{dt}\int|\mathbf{b}|^2\bar{x}^adx
+\nu\int|\nabla\mathbf{b}|^2\bar{x}^adx
& = \frac{\nu}{2}\int|\mathbf{b}|^2\Delta\bar{x}^adx
+\int(\mathbf{b}\cdot\nabla)\mathbf{u}\cdot\mathbf{b}\bar{x}^adx
\notag \\
& \quad -\frac12\int\divv\mathbf{u}|\mathbf{b}|^2\bar{x}^adx
+\frac12\int|\mathbf{b}|^2\mathbf{u}\cdot\nabla\bar{x}^adx
\notag \\
& \triangleq\sum_{i=1}^4K_i.
\end{align}
Direct calculations lead to
\begin{equation}\label{10.4}
|K_1|\leq C\int|\mathbf{b}|^2\bar{x}^a\bar{x}^{-2}\log^{2(1-\eta_0)}(e+|x|^2)dx
\leq C\|\mathbf{b}\bar{x}^{\frac{a}{2}}\|_{L^2}^2,
\end{equation}
and
\begin{align}\label{10.5}
|K_2|+|K_3| & \leq C\int|\nabla\mathbf{u}||\mathbf{b}|^2\bar{x}^adx \notag \\
&\leq C\|\nabla\mathbf{u}\|_{L^2}\|\mathbf{b}\bar{x}^{\frac{a}{2}}\|_{L^4}^2
 \notag \\
&\leq C\|\nabla\mathbf{u}\|_{L^2}\|\mathbf{b}\bar{x}^{\frac{a}{2}}\|_{L^2}
\left(\|\nabla\mathbf{b}\bar{x}^{\frac{a}{2}}\|_{L^2}
+\|\mathbf{b}\nabla\bar{x}^{\frac{a}{2}}\|_{L^2}\right)
 \notag \\
&\leq C\|\mathbf{b}\bar{x}^{\frac{a}{2}}\|_{L^2}
\left(\|\nabla\mathbf{b}\bar{x}^{\frac{a}{2}}\|_{L^2}
+\|\mathbf{b}\bar{x}^{\frac{a}{2}}\|_{L^2}
\|\bar{x}^{-1}\nabla\bar{x}\|_{L^\infty}\right)
 \notag \\
&\leq C\|\mathbf{b}\bar{x}^{\frac{a}{2}}\|_{L^2}^2
+\frac{\nu}{4}\|\nabla\mathbf{b}\bar{x}^{\frac{a}{2}}\|_{L^2}^2,
\end{align}
due to
\begin{align*}
|\nabla\bar{x}|\leq (3+2\eta_0)\log^{1+\eta_0}(e+|x|^2)
\leq C(a,\eta_0)\bar{x}^{\frac{4}{8+a}}.
\end{align*}
It follows from H{\"o}lder's inequality, Gagliardo-Nirenberg inequality, and \eqref{06.2} that
\begin{align}\label{10.6}
|K_4| & \leq C\int|\mathbf{b}|^2\bar{x}^a
\bar{x}^{-\frac34}|\mathbf{u}|\bar{x}^{-\frac14}\log^{(1-\eta_0)}(e+|x|^2)dx \notag \\
&\leq C\|\mathbf{b}\bar{x}^{\frac{a}{2}}\|_{L^4}
\|\mathbf{b}\bar{x}^{\frac{a}{2}}\|_{L^2}
\|\mathbf{u}\bar{x}^{-\frac{3}{4}}\|_{L^4}
\notag \\
&\leq C\|\mathbf{b}\bar{x}^{\frac{a}{2}}\|_{L^2}^2
+\frac{\nu}{4}\|\nabla\mathbf{b}\bar{x}^{\frac{a}{2}}\|_{L^2}^2.
\end{align}
Putting \eqref{10.4}--\eqref{10.6} into \eqref{10.3} and using Gronwall's inequality, we obtain the desired \eqref{10.1}.

2. Multiplying \eqref{1.10}$_4$ by $\Delta\mathbf{b}\bar{x}^a$ and integrating the resulting equations yield that
\begin{align}\label{10.7}
& \frac12\frac{d}{dt}\int|\nabla\mathbf{b}|^2\bar{x}^adx
+\nu\int|\Delta\mathbf{b}|^2\bar{x}^adx \notag \\
& \leq C\int|\nabla\mathbf{b}||\mathbf{b}||\nabla\mathbf{u}||\nabla\bar{x}^a|dx
+C\int|\nabla\mathbf{b}|^2|\mathbf{u}||\nabla\bar{x}^a|dx
+C\int|\nabla\mathbf{b}||\Delta\mathbf{b}||\nabla\bar{x}^a|dx
\notag \\
& \quad +C\int|\mathbf{b}||\nabla\mathbf{u}||\Delta\mathbf{b}|\bar{x}^adx
+C\int|\nabla\mathbf{u}||\nabla\mathbf{b}|^2\bar{x}^adx
\triangleq\sum_{i=1}^5J_i.
\end{align}
Applying Gagliardo-Nirenberg inequality and \eqref{10.1}, we have
\begin{align}
|J_1| & \leq C\int|\nabla\mathbf{b}||\mathbf{b}|\nabla\mathbf{u}|
\bar{x}^a(\bar{x}^{-1}|\nabla\bar{x}|)dx \notag \\
&\leq C\|\mathbf{b}\bar{x}^{\frac{a}{2}}\|_{L^4}^4
+C\|\nabla\mathbf{b}\bar{x}^{\frac{a}{2}}\|_{L^2}^2
+C\|\nabla\mathbf{u}\|_{L^4}^4
\notag \\
&\leq C\|\mathbf{b}\bar{x}^{\frac{a}{2}}\|_{L^2}^2
\left(\|\nabla\mathbf{b}\bar{x}^{\frac{a}{2}}\|_{L^2}^2
+\|\mathbf{b}\bar{x}^{\frac{a}{2}}\|_{L^2}^2\right)
+C\|\nabla\mathbf{b}\bar{x}^{\frac{a}{2}}\|_{L^2}^2
+C\|\nabla\mathbf{u}\|_{L^4}^4\notag \\
&\leq C\|\nabla\mathbf{b}\bar{x}^{\frac{a}{2}}\|_{L^2}^2
+C\|\nabla\mathbf{u}\|_{L^4}^4;\label{10.8} \\
|J_2| & \leq C\int|\nabla\mathbf{b}|^{\frac{4a-1}{2a}}
\bar{x}^{\frac{4a-1}{4}}|\nabla\mathbf{b}|^{\frac{1}{2a}}
|\mathbf{u}|\bar{x}^{-\frac12}
\bar{x}^{-\frac14}|\nabla\bar{x}|dx \notag \\
& \leq C\||\nabla\mathbf{b}|^{\frac{4a-1}{2a}}
\bar{x}^{\frac{4a-1}{4}}\|_{L^{\frac{4a}{4a-1}}}
\||\nabla\mathbf{b}|^{\frac{1}{2a}}\|_{L^{8a}}
\|\mathbf{u}\bar{x}^{-\frac12}\|_{L^{8a}} \notag \\
&\leq C\|\nabla\mathbf{b}\bar{x}^{\frac{a}{2}}\|_{L^2}^2
+C\|\nabla\mathbf{b}\|_{L^4}^2
\notag \\
&\leq C\|\nabla\mathbf{b}\bar{x}^{\frac{a}{2}}\|_{L^2}^2
+C\|\nabla^2\mathbf{b}\|_{L^2}^2
+C; \label{10.9} \\
|J_3|+|J_4| &\leq C\|\Delta\mathbf{b}\bar{x}^{\frac{a}{2}}\|_{L^2}
\|\nabla\mathbf{b}\bar{x}^{\frac{a}{2}}\|_{L^2}
+C\|\Delta\mathbf{b}\bar{x}^{\frac{a}{2}}\|_{L^2}
\|\mathbf{b}\bar{x}^{\frac{a}{2}}\|_{L^4}\|\nabla\mathbf{u}\|_{L^4}
\notag \\
& \leq \frac{\nu}{4}\|\Delta\mathbf{b}\bar{x}^{\frac{a}{2}}\|_{L^2}^2
+C\|\nabla\mathbf{b}\bar{x}^{\frac{a}{2}}\|_{L^2}^2
+C\|\nabla\mathbf{u}\|_{L^4}^4; \label{10.10} \\
|J_5| & \leq C\|\nabla\mathbf{u}\|_{L^2}
\|\nabla\mathbf{b}\bar{x}^{\frac{a}{2}}\|_{L^4}^2 \notag \\
& \leq C\|\nabla\mathbf{b}\bar{x}^{\frac{a}{2}}\|_{L^2}
\left(\|\nabla^2\mathbf{b}\bar{x}^{\frac{a}{2}}\|_{L^2}
+\|\nabla\mathbf{b}\bar{x}^{\frac{a}{2}}\|_{L^2}
\|\bar{x}^{-1}\nabla\bar{x}\|_{L^\infty}\right)\notag \\
& \leq \frac{\nu}{4}\|\nabla^2\mathbf{b}\bar{x}^{\frac{a}{2}}\|_{L^2}^2
+C\|\nabla\mathbf{b}\bar{x}^{\frac{a}{2}}\|_{L^2}^2.
\label{10.11}
\end{align}
Inserting \eqref{10.8}--\eqref{10.11} into \eqref{10.7},
and noting the following fact
\begin{align*}
\int|\nabla^2\mathbf{b}|^2\bar{x}^adx
& = \int|\Delta\mathbf{b}|^2\bar{x}^adx
-\int\partial_i\partial_k\mathbf{b}\cdot\partial_k\mathbf{b}\partial_i\bar{x}^adx
+\int\partial_i\partial_i\mathbf{b}\cdot\partial_k\mathbf{b}\partial_k\bar{x}^adx \\
& \leq \int|\Delta\mathbf{b}|^2\bar{x}^adx+\frac12\int|\nabla^2\mathbf{b}|^2\bar{x}^adx
+C\int|\nabla\mathbf{b}|^2\bar{x}^adx,
\end{align*}
we derive that
\begin{align}\label{10.12}
\frac{d}{dt}\int|\nabla\mathbf{b}|^2\bar{x}^adx
+\int|\nabla^2\mathbf{b}|^2\bar{x}^adx
\leq C\|\nabla\mathbf{b}\bar{x}^{\frac{a}{2}}\|_{L^2}^2
+C\|\nabla^2\mathbf{b}\|_{L^2}^2+C\|\nabla\mathbf{u}\|_{L^4}^4.
\end{align}
It follows from \eqref{2.5} and \eqref{3.1} that
\begin{align}\label{6.17}
\|\nabla\mathbf{u}\|_{L^4}^{4}
\leq C(\|\rho\dot{\mathbf{u}}\|_{L^2}+\||\mathbf{b}||\nabla\mathbf{b}|\|_{L^2})^2
+C \leq C\|\sqrt{\rho}\dot{\mathbf{u}}\|_{L^2}^{2}
+C\||\mathbf{b}||\nabla\mathbf{b}|\|_{L^2}^2+C,
\end{align}
which together with \eqref{10.12} leads to
\begin{align*}
\frac{d}{dt}\|\nabla\mathbf{b}\bar{x}^{\frac{a}{2}}\|_{L^2}^2
+\|\nabla^2\mathbf{b}\bar{x}^{\frac{a}{2}}\|_{L^2}^2
\leq C\|\nabla\mathbf{b}\bar{x}^{\frac{a}{2}}\|_{L^2}^2
+C\|\nabla^2\mathbf{b}\|_{L^2}^2+C\|\sqrt{\rho}\dot{\mathbf{u}}\|_{L^2}^{2}
+C\||\mathbf{b}||\nabla\mathbf{b}|\|_{L^2}^2+C.
\end{align*}
This along with Gronwall's inequality, \eqref{5.1},  and \eqref{3.5} yields the desired \eqref{10.2}.
\hfill $\Box$

\begin{lemma}\label{lem35}
Under the condition \eqref{3.1}, it holds that for any $T\in[0,T^*)$,
\begin{equation}\label{6.1}
\sup_{0\leq t\leq T}\left(\|\mathbf{b}_t\|_{L^2}^2+\|\nabla^2\mathbf{b}\|_{L^2}^2
+\|\sqrt{\rho}\dot{\mathbf{u}}\|_{L^2}^2\right)
+\int_0^T\left(\|\nabla\mathbf{b}_t\|_{L^2}^2
+\|\nabla\dot{\mathbf{u}}\|_{L^2}^2\right)dt \leq C.
\end{equation}
\end{lemma}
{\it Proof.}
1. Differentiating $\eqref{1.10}_4$ with respect to $t$, we have
\begin{equation}\label{11.2}
\mathbf{b}_{tt}-\nu\Delta\mathbf{b}_t
=\mathbf{b}_t\cdot\nabla\mathbf{u}
-\mathbf{u}\cdot\nabla\mathbf{b}_t
-\mathbf{b}_t\divv\mathbf{u}+\mathbf{b}\cdot\nabla\mathbf{u}_t
-\mathbf{u}_t\cdot\nabla\mathbf{b}-\mathbf{b}\divv\mathbf{u}_t.
\end{equation}
Multiplying \eqref{11.2} by $\mathbf{b}_t$ and integrating by parts lead to
\begin{align}\label{11.3}
\frac12\frac{d}{dt}\int|\mathbf{b}_t|^2dx+\nu\int|\nabla\mathbf{b}_t|^2dx
&=\int(\mathbf{b}_t\cdot\nabla\mathbf{u}-\mathbf{u}\cdot\nabla\mathbf{b}_t
-\mathbf{b}_t\divv\mathbf{u})\cdot\mathbf{b}_tdx \nonumber\\
& \quad +\int(\mathbf{b}\cdot\nabla\mathbf{u}_t-\mathbf{u}_t\cdot\nabla\mathbf{b}
-\mathbf{b}\divv\mathbf{u}_t)\cdot\mathbf{b}_tdx\nonumber\\
&\triangleq L_1+L_2.
\end{align}
It follows from integration by parts, Gagliardo-Nirenberg inequality,  and \eqref{5.1} that
\begin{align}\label{11.4}
L_1 =\int \left(\mathbf{b}_t\cdot\nabla\mathbf{u}\cdot\mathbf{b}_t
-\frac12|\mathbf{b}_t|^2\divv\mathbf{u}\right)dx
\leq C\|\nabla\mathbf{u}\|_{L^2}\|\mathbf{b}_t\|_{L^4}^2
\leq \frac{\nu}{4}\|\nabla\mathbf{b}_t\|_{L^2}^2
+C\|\mathbf{b}_t\|_{L^2}^2.
\end{align}
Similarly to \eqref{06.2}, one infers from \eqref{22}, \eqref{p1}, and \eqref{3.1} that for any $\eta\in(0,1]$ and any $s>2$,
\begin{equation}\label{06.3}
\|\dot{\mathbf{u}}\bar{x}^{-\eta}\|_{L^{s/\eta}}\leq C\left(\|\sqrt{\rho}\dot{\mathbf{u}}\|_{L^2}
+\|\nabla\dot{\mathbf{u}}\|_{L^2}\right),
\end{equation}
which combined with H{\"o}lder's inequality, Gagliardo-Nirenberg inequality, \eqref{3.5}, \eqref{10.1}, \eqref{10.2}, and \eqref{6.17} leads to
\begin{align}\label{11.5}
L_2 & =\int (\mathbf{b}\cdot\nabla\dot{\mathbf{u}}-\dot{\mathbf{u}}\cdot\nabla\mathbf{b}
-\mathbf{b}\divv\dot{\mathbf{u}})\cdot\mathbf{b}_tdx   \nonumber\\
& \quad -\int [\mathbf{b}\cdot\nabla(\mathbf{u}\cdot\nabla\mathbf{u})
-(\mathbf{u}\cdot\nabla\mathbf{u})\cdot\nabla\mathbf{b}
-\mathbf{b}\divv(\mathbf{u}\cdot\nabla\mathbf{u})]\cdot\mathbf{b}_tdx  \nonumber\\
&= \int (\mathbf{b}\cdot\nabla\dot{\mathbf{u}}-\dot{\mathbf{u}}\cdot\nabla\mathbf{b}
-\mathbf{b}\divv\dot{\mathbf{u}})\cdot\mathbf{b}_t\mbox{d}x   \nonumber\\
& \quad +\int [(\mathbf{u}\cdot\nabla\mathbf{u})\cdot(\mathbf{b}\cdot\nabla\mathbf{b}_t)
+(\mathbf{u}\cdot\nabla\mathbf{u})\cdot\nabla\mathbf{b}_t\cdot\mathbf{b}]dx   \nonumber\\
&\leq C\int \left(|\mathbf{b}||\mathbf{b}_t||\nabla\dot{\mathbf{u}}|
+|\dot{\mathbf{u}}||\nabla\mathbf{b}||\mathbf{b}_t|
+|\mathbf{u}||\nabla\mathbf{u}||\mathbf{b}||\nabla\mathbf{b}_t|\right)dx  \nonumber\\
&\leq C\|\mathbf{b}\|_{L^4}\|\mathbf{b}_t\|_{L^4}\|\nabla\dot{\mathbf{u}}\|_{L^2}
+C\|\dot{\mathbf{u}}\bar{x}^{-\frac{a}{2}}\|_{L^2}
\|\nabla\mathbf{b}\bar{x}^{\frac{a}{2}}\|_{L^4}\|\mathbf{b}_t\|_{L^4}\nonumber\\
& \quad
+C\|\mathbf{u}\bar{x}^{-\frac{a}{2}}\|_{L^8}
\|\nabla\mathbf{u}\|_{L^4}\|\mathbf{b}\bar{x}^{\frac{a}{2}}\|_{L^8}
\|\nabla\mathbf{b}_t\|_{L^2}  \nonumber\\
& \leq \frac{\nu}{4}\|\nabla\mathbf{b}_t\|_{L^2}^2
+C\|\nabla\dot{\mathbf{u}}\|_{L^2}^2
+C\|\mathbf{b}_t\|_{L^2}^2+C\|\nabla^2\mathbf{b}\bar{x}^{\frac{a}{2}}\|_{L^2}^2
+C\|\sqrt{\rho}\dot{\mathbf{u}}\|_{L^2}^{2}
+C\||\mathbf{b}||\nabla\mathbf{b}|\|_{L^2}^2+C.
\end{align}
Inserting \eqref{11.4} and \eqref{11.5} into \eqref{11.3}, we have
\begin{align}\label{eq11.6}
\frac{d}{dt}\|\mathbf{b}_t\|_{L^2}^2+\nu\|\nabla\mathbf{b}_t\|_{L^2}^2
& \leq C\|\mathbf{b}_t\|_{L^2}^2
+C_1\|\nabla\dot{\mathbf{u}}\|_{L^2}^2
+C\|\nabla^2\mathbf{b}\bar{x}^{\frac{a}{2}}\|_{L^2}^2 \notag \\
& \quad
+C\|\sqrt{\rho}\dot{\mathbf{u}}\|_{L^2}^{2}
+C\||\mathbf{b}||\nabla\mathbf{b}|\|_{L^2}^2+C.
\end{align}

2. Operating $\partial_t+\divv(\mathbf{u}\cdot)$ to $j$-th component of \eqref{1.10}$_2$
and multiplying the resulting equation by $\dot{u}_j$, one gets by some calculations that
\begin{align}\label{6.2}
\frac12\frac{d}{dt}\int \rho|\dot{\mathbf{u}}|^2\mbox{d}x
& =\mu\int\dot{u}_j(\partial_t\Delta u_j+\divv(\mathbf{u}\Delta u_j))dx
+(\lambda+\mu)\int \dot{u}_j (\partial_t\partial_j(\divv\mathbf{u})+\divv(\mathbf{u}\partial_j(\divv\mathbf{u})))dx  \nonumber\\
& \quad -\int \dot{u}_j(\partial_jP_t+\divv(\mathbf{u}\partial_jP)) dx
-\frac12\int \dot{u}_j (\partial_t\partial_j|\mathbf{b}|^2+\divv(\mathbf{u}\partial_j|\mathbf{b}|^2)) dx\nonumber\\
& \quad +\int \dot{u}_j (\partial_t(\mathbf{b}\cdot\nabla b_j)+\divv(\mathbf{u}(\mathbf{b}\nabla b_j))) dx
\triangleq\sum_{i=1}^{5}J_i.
\end{align}
Integration by parts leads to
\begin{align}\label{6.3}
J_1 & = - \mu\int(\partial_i\dot{u}_j\partial_t\partial_i u_j+\Delta u_j\mathbf{u}\cdot\nabla\dot{u}_j)dx \notag \\
& = - \mu\int(|\nabla\dot{\mathbf{u}}|^2
-\partial_i\dot{u}_ju_k\partial_k\partial_i u_j
-\partial_i\dot{u}_j\partial_iu_k\partial_k u_j
+\Delta u_j\mathbf{u}\cdot\nabla\dot{u}_j)dx
\notag \\
& = - \mu\int(|\nabla\dot{\mathbf{u}}|^2
+\partial_i\dot{u}_j\partial_k u_k\partial_i u_j
-\partial_i\dot{u}_j\partial_iu_k\partial_k u_j
-\partial_i u_j\partial_iu_k\partial_k\dot{u}_j)dx
\notag \\
& \leq -\frac{3\mu}{4}\|\nabla\dot{\mathbf{u}}\|_{L^2}^2+C\|\nabla\mathbf{u}\|_{L^4}^4.
\end{align}
Similarly, one has
\begin{align}\label{6.4}
J_2 \leq -\frac{\lambda+\mu}{2}\|\divv\dot{\mathbf{u}}\|_{L^2}^2
+C\|\nabla\mathbf{u}\|_{L^4}^4.
\end{align}
It follows from integration by parts, \eqref{1.10}$_3$, \eqref{3.1}, and \eqref{5.1} that
\begin{align}\label{6.5}
J_3 & = \int(\partial_j\dot{u}_jP_t+\partial_j P\mathbf{u}\cdot\nabla\dot{u}_j)dx \notag \\
& = \int\partial_j\dot{u}_j\left(
2\mu|\mathfrak{D}(\mathbf{u})|^2+\lambda(\divv\mathbf{u})^2
+\nu|\nabla\times\mathbf{b}|^2-\divv(P\mathbf{u})
-P\divv\mathbf{u}\right)dx \notag \\
& \quad -\int P\partial_j(\mathbf{u}\cdot\nabla\dot{u}_j)dx\notag \\
& = \int\partial_j\dot{u}_j\left(
2\mu|\mathfrak{D}(\mathbf{u})|^2+\lambda(\divv\mathbf{u})^2
+\nu|\nabla\times\mathbf{b}|^2-\divv(P\mathbf{u})
-P\divv\mathbf{u}\right)dx \notag \\
& \quad
-\int P\partial_j\mathbf{u}\cdot\nabla\dot{u}_jdx
+\int\partial_j\dot{u}_j \divv(P\mathbf{u})dx
\notag \\
& = \int\partial_j\dot{u}_j\left(
2\mu|\mathfrak{D}(\mathbf{u})|^2+\lambda(\divv\mathbf{u})^2
+\nu|\nabla\times\mathbf{b}|^2
-P\divv\mathbf{u}\right)dx-\int P\partial_j\mathbf{u}\cdot\nabla\dot{u}_jdx
\notag \\
& \leq C\int |\nabla\dot{\mathbf{u}}|(|\nabla\mathbf{u}|^2+|\nabla\mathbf{b}|^2+1)dx
\notag \\
& \leq \frac{\mu}{4}\|\nabla\dot{\mathbf{u}}\|_{L^2}^2
+C\|\nabla\mathbf{u}\|_{L^4}^4+C\|\nabla\mathbf{b}\|_{L^4}^4+C.
\end{align}
From \eqref{1.10}$_4$, \eqref{1.10}$_5$, \eqref{3.5}, and Gagliardo-Nirenberg inequality, we arrive at
\begin{align}\label{6.6}
J_4 & = \int\partial_j\dot{u}_j\mathbf{b}\cdot\mathbf{b}_tdx +\frac12\int\mathbf{u}\cdot\nabla\dot{u}_j\partial_j|\mathbf{b}|^2dx \notag \\
& = \int\partial_j\dot{u}_j\mathbf{b}\cdot\mathbf{b}_tdx +\frac12\int\partial_j\dot{u}_j\divv\mathbf{u}|\mathbf{b}|^2dx
-\frac12\int\partial_j\mathbf{u}\cdot\nabla\dot{u}_j|\mathbf{b}|^2dx \notag \\
& \leq C\int |\nabla\dot{\mathbf{u}}||\mathbf{b}||\mathbf{b}_t|dx
+C\int |\nabla\dot{\mathbf{u}}||\nabla\mathbf{u}||\mathbf{b}|^2dx
\notag \\
& \leq \frac{\mu}{8}\|\nabla\dot{\mathbf{u}}\|_{L^2}^2
+C\|\nabla\mathbf{u}\|_{L^4}^4+C\|\mathbf{b}_t\|_{L^2}^2
+\varepsilon\|\nabla\mathbf{b}_t\|_{L^2}^2.
\end{align}
Similarly to $J_4$, we infer that
\begin{align}\label{6.7}
J_5 \leq \frac{\mu}{8}\|\nabla\dot{\mathbf{u}}\|_{L^2}^2
+C\|\nabla\mathbf{u}\|_{L^4}^4+C\|\mathbf{b}_t\|_{L^2}^2
+\varepsilon\|\nabla\mathbf{b}_t\|_{L^2}^2.
\end{align}
Inserting \eqref{6.3}--\eqref{6.7} into \eqref{6.2} and applying Garliardo-Nirenberg inequality and \eqref{6.17}, one has
\begin{align}\label{6.8}
\frac{d}{dt}\|\sqrt{\rho}\dot{\mathbf{u}}\|_{L^2}^2
+\mu\|\nabla\dot{\mathbf{u}}\|_{L^2}^2
& \leq C\|\nabla\mathbf{u}\|_{L^4}^4+C\|\nabla\mathbf{b}\|_{L^4}^4
+C\|\mathbf{b}_t\|_{L^2}^2
+C\|\nabla\mathbf{b}_t\|_{L^2}^2 \notag \\
& \leq C\|\sqrt{\rho}\dot{\mathbf{u}}\|_{L^2}^{2}
+C\||\mathbf{b}||\nabla\mathbf{b}|\|_{L^2}^2
+C\|\nabla^2\mathbf{b}\|_{L^2}^2+C\|\mathbf{b}_t\|_{L^2}^2 \notag \\
& \quad
+2\varepsilon\|\nabla\mathbf{b}_t\|_{L^2}^2+C.
\end{align}
Then, adding \eqref{6.8} multiplied by $\frac{C_1+1}{\mu}$ to \eqref{eq11.6} and choosing $\varepsilon$ suitably small, we derive that
\begin{align*}
& \frac{d}{dt}\left(\frac{C_1+1}{\mu}\|\sqrt{\rho}\dot{\mathbf{u}}\|_{L^2}^2
+\|\mathbf{b}_t\|_{L^2}^2\right)
+\|\nabla\dot{\mathbf{u}}\|_{L^2}^2+\nu\|\nabla\mathbf{b}_t\|_{L^2}^2 \notag \\
& \leq C\left(\|\sqrt{\rho}\dot{\mathbf{u}}\|_{L^2}^{2}
+\|\mathbf{b}_t\|_{L^2}^2\right)
+C\||\mathbf{b}||\nabla\mathbf{b}|\|_{L^2}^2
+C\|\nabla^2\mathbf{b}\|_{L^2}^2
+C\|\nabla^2\mathbf{b}\bar{x}^{\frac{a}{2}}\|_{L^2}^2+C.
\end{align*}
This together with Gronwall's inequality, \eqref{3.5}, \eqref{5.1},  and \eqref{10.2} implies
\begin{equation}\label{eq6.1}
\sup_{0\leq t\leq T}\left(\|\sqrt{\rho}\dot{\mathbf{u}}\|_{L^2}^2
+\|\mathbf{b}_t\|_{L^2}^2\right)
+\int_0^T\left(\|\nabla\dot{\mathbf{u}}\|_{L^2}^2
+\|\nabla\mathbf{b}_t\|_{L^2}^2\right)dt \leq C.
\end{equation}

3. Applying the standard $L^2$-estimate to \eqref{1.10}$_4$ and using Gagliardo-Nirenberg inequality, \eqref{10.2}, \eqref{5.1}, and \eqref{3.5}, we get
\begin{align*}
\|\nabla^2 \mathbf{b}\|_{L^2} & \leq C\left(\|\mathbf{b}_t\|_{L^2}
+\||\mathbf{u}||\nabla\mathbf{b}|\|_{L^2}
+\||\mathbf{b}||\nabla\mathbf{u}|\|_{L^2}\right)
\nonumber\\
&\leq C\left(\|\mathbf{b}_t\|_{L^2}
+\|\mathbf{u}\bar{x}^{-\frac{a}{4}}\|_{L^8}^{2}
\||\nabla\mathbf{b}|^{\frac12}\bar{x}^{\frac{a}{4}}\|_{L^4}^2
\||\nabla\mathbf{b}|^{\frac12}\|_{L^8}^2
+\|\mathbf{b}\|_{L^\infty}\|\nabla\mathbf{u}\|_{L^2}\right)\nonumber\\
&\leq C+C\|\nabla\mathbf{b}\bar{x}^{\frac{a}{2}}\|_{L^2}
\|\nabla\mathbf{b}\|_{L^4}+C\|\mathbf{b}\|_{L^4}^{\frac12}
\|\nabla\mathbf{b}\|_{L^4}^{\frac12} \nonumber\\
&\leq \frac12
\|\nabla^2\mathbf{b}\|_{L^2}+C.
\end{align*}
Thus, one has
\begin{equation}\label{11.7}
\sup_{0\leq t\leq T}\|\nabla^2\mathbf{b}\|_{L^2}^2\leq C.
\end{equation}
The proof of Lemma \ref{lem35} is completed.   \hfill $\Box$

The following lemma will treat the higher order derivatives of the solutions which are needed to guarantee the extension of local strong solution to be a global one.
\begin{lemma}\label{lem37}
Under the condition \eqref{3.1}, and let $q>2$ be as in Theorem \ref{thm1.1}, then it holds that for any $T\in[0,T^*)$,
\begin{equation}\label{7.1}
\sup_{0\leq t\leq T}\left(\|(\rho,P)\|_{H^1\cap W^{1,q}}+\|\nabla\mathbf{u}\|_{H^1}\right)\leq C.
\end{equation}
\end{lemma}
{\it Proof.}
1. It follows from the mass equation \eqref{1.10}$_1$ that $\nabla\rho$ satisfies for any $r\in[2,q]$,
\begin{align}\label{7.2}
\frac{d}{dt}\|\nabla\rho\|_{L^r}
& \leq C(r)(1+\|\nabla\mathbf{u} \|_{L^\infty})\|\nabla\rho\|_{L^r}
+C(r)\|\nabla^2\mathbf{u}\|_{L^r} \nonumber \\
& \leq C(1+\|\nabla\mathbf{u} \|_{L^\infty})\|\nabla\rho\|_{L^r}
+C(\|\rho\dot{\mathbf{u}}\|_{L^r}+\|\nabla P\|_{L^r}+\||\mathbf{b}||\nabla\mathbf{b}|\|_{L^r})
\end{align}
due to
\begin{equation}\label{7.3}
\|\nabla^2\mathbf{u}\|_{L^r}
\leq C(\|\rho\dot{\mathbf{u}}\|_{L^r}+\|\nabla P\|_{L^r}+\||\mathbf{b}||\nabla\mathbf{b}|\|_{L^r}),
\end{equation}
which follows from the standard $L^r$-estimate for the following elliptic system
\begin{align*}
\begin{cases}
\mu\Delta\mathbf{u}+(\lambda+\mu)\nabla\divv\mathbf{u}=\rho\dot{\mathbf{u}}+\nabla P-\mathbf{b}\cdot\nabla\mathbf{b}+\frac12\nabla|\mathbf{b}|^2,\ \ x\in\mathbb{R}^2,\\
\mathbf{u}\rightarrow\mathbf{0},\ \text{as}\ |x|\rightarrow\infty.
\end{cases}
\end{align*}
Similarly, one deduces from \eqref{1.10}$_3$ that $\nabla P$ satisfies for any $r\in[2,q]$,
\begin{align}\label{7.4}
\frac{d}{dt}\|\nabla P\|_{L^r} & \leq C(r)(1+\|\nabla\mathbf{u}\|_{L^{\infty}})(\|\nabla P\|_{L^r}+\|\nabla^2\mathbf{u}\|_{L^r})
+C(r)\|\nabla\mathbf{b}\|_{L^{\infty}}\|\nabla^2\mathbf{b}\|_{L^r} \notag \\
& \leq C(1+\|\nabla\mathbf{u}\|_{L^{\infty}})(\|\rho\dot{\mathbf{u}}\|_{L^r}+\|\nabla P\|_{L^r}+\||\mathbf{b}||\nabla\mathbf{b}|\|_{L^r})
+C\|\nabla\mathbf{b}\|_{L^\infty}
\|\nabla^2\mathbf{b}\|_{L^r}.
\end{align}

2. We infer from Sobolev's inequality, \eqref{3.3}, \eqref{5.1}, and \eqref{11.7} that
\begin{equation}\label{eq7.06}
\sup_{0\leq t\leq T}\|\mathbf{b}\|_{L^\infty}
\leq\sup_{0\leq t\leq T}\|\mathbf{b}\|_{H^2}\leq C,
\end{equation}
which combined with \eqref{3.14}, Gagliardo-Nirenberg inequality, \eqref{3.1}, and \eqref{2.3} yields
\begin{align}\label{7.5}
\|\divv\mathbf{u}\|_{L^\infty}+\|\omega\|_{L^\infty}
& \leq C\|P\|_{L^\infty}+C\|F\|_{L^\infty}
+C\||\mathbf{b}|^2\|_{L^\infty}+\|\omega\|_{L^\infty} \nonumber \\
& \leq C+C(q)\|\nabla F\|_{L^2}^{\frac{q-2}{2(q-1)}}\|\nabla F\|_{L^q}^{\frac{q}{2(q-1)}}
+C(q)\|\nabla\omega\|_{L^2}^{\frac{q-2}{2(q-1)}}\|\nabla\omega\|_{L^q}^{\frac{q}{2(q-1)}} \nonumber \\
& \leq C+C\|\rho\dot{\mathbf{u}}\|_{L^q}^{\frac{q}{2(q-1)}}
+C\||\mathbf{b}||\nabla\mathbf{b}|\|_{L^q}^{\frac{q}{2(q-1)}}
\nonumber \\
& \leq C+C\|\rho\dot{\mathbf{u}}\|_{L^q}^{\frac{q}{2(q-1)}}
+C\|\mathbf{b}\|_{L^\infty}^{\frac{q}{2(q-1)}}
\|\nabla\mathbf{b}\|_{H^1}^{\frac{q}{2(q-1)}}\nonumber \\
& \leq C+C\|\rho\dot{\mathbf{u}}\|_{L^q}^{\frac{q}{2(q-1)}}.
\end{align}
This along with Lemma \ref{lem25}, \eqref{7.3}, and \eqref{5.1} gives rise to
\begin{align}\label{7.6}
\|\nabla\mathbf{u}\|_{L^\infty}
& \leq C\left(\|\divv\mathbf{u}\|_{L^\infty}+\|\omega\|_{L^\infty}\right)
\log(e+\|\nabla^2\mathbf{u}\|_{L^q})+C\|\nabla\mathbf{u}\|_{L^2}+C \nonumber \\
& \leq C\left(1+\|\rho\dot{\mathbf{u}}\|_{L^q}^{\frac{q}{2(q-1)}}\right)
\log\left(e+\|\rho\dot{\mathbf{u}}\|_{L^q}+\|\nabla P\|_{L^q}\right)+C.
\end{align}
It follows from \eqref{p1}, \eqref{3.1}, \eqref{22}, and \eqref{06.1} that for any $\eta\in(0,1]$ and any $s>2$,
\begin{align}\label{7.9}
\|\rho^\eta v\|_{L^{\frac{s}{\eta}}}
 &  \leq C\|\rho^\eta\bar x^{\frac{3\eta a}{4s}}\|_{L^{\frac{4s}{3\eta}}}
\|v\bar  x^{-\frac{3\eta a}{4s}}\|_{L^{\frac{4s}{\eta}}} \notag \\
 &  \leq C\|\rho\|_{L^\infty}^{\frac{(4s-3)\eta}{4s}}\|\rho\bar x^a\|_{L^1}^{\frac{3\eta}{4s}}\left( \|\sqrt{\rho} v\|_{L^2}+\|\nabla v\|_{L^2}\right) \notag \\
 &  \leq C\left(\|\sqrt{\rho}v\|_{L^2}+\|\nabla v\|_{L^2}\right),
\end{align}
which together with H{\"o}lder's inequality, \eqref{6.1}, and \eqref{3.1} shows that
\begin{align}\label{7.10}
\|\rho\dot{\mathbf{u}}\|_{L^q}
 &  \leq C\|\rho\dot{\mathbf{u}}\|_{L^2}^{\frac{2(q-1)}{q^2-2}}
\|\rho\dot{\mathbf{u}}\|_{L^{q^2}}^{\frac{q(q-2)}{q^2-2}} \notag \\
 &  \leq C\|\rho\dot{\mathbf{u}}\|_{L^2}^{\frac{2(q-1)}{q^2-2}}
\left(\|\sqrt{\rho}\dot{\mathbf{u}}\|_{L^2}+\|\nabla\dot{\mathbf{u} }\|_{L^2}\right)^{\frac{q(q-2)}{q^2-2}} \notag \\
 &  \leq C\left(1+\|\nabla\dot{\mathbf{u}}\|_{L^2}^{\frac{q(q-2)}{q^2-2}}\right),
\end{align}
Then we derive from \eqref{7.6} and \eqref{7.10} that
\begin{align}\label{7.11}
\|\nabla\mathbf{u}\|_{L^\infty}
\leq C\left(1+\|\nabla\dot{\mathbf{u}}\|_{L^2}\right)
\log\left(e+\|\nabla\dot{\mathbf{u}}\|_{L^2}+\|\nabla P\|_{L^q}\right)+C
\end{align}
due to $\frac{q(q^2-2q)}{(2q-2)(q^2-2)},\ \frac{q^2-2q}{q^2-2}\in(0,1)$.
Hence, substituting \eqref{7.10} and \eqref{7.11} into \eqref{7.2} and \eqref{7.4}, we get after choosing $r=q$ that
\begin{equation}\label{7.011}
f'(t)\leq  Cg(t)f(t)\log{f(t)}+Cg(t)f(t)+Cg(t),
\end{equation}
where
\begin{align*}
f(t)&\triangleq e+\|\nabla\rho\|_{L^{ q}}+\|\nabla P\|_{L^{ q}},\\
g(t)&\triangleq (1+\|\nabla\dot{\mathbf{u}}\|_{L^2})
\log(e+\|\nabla\dot{\mathbf{u}}\|_{L^2})+\|\nabla\mathbf{b}_t\|_{L^2}^2.
\end{align*}
This yields
\begin{equation}\label{7.12}
(\log f(t))'\leq Cg(t)+Cg(t)\log f(t)
\end{equation}
due to $f(t)>1$.
Thus it follows from \eqref{7.12}, \eqref{6.1}, and Gronwall's inequality that
\begin{equation}\label{7.13}
\sup_{0\leq t\leq T}\left(\|\nabla \rho\|_{L^q}+\|\nabla P\|_{L^q}\right)\leq C.
\end{equation}

3. Taking $r=2$ in \eqref{7.2} and \eqref{7.4}, one gets from \eqref{3.1}, \eqref{6.1}, \eqref{eq7.06},  and Garliardo-Nirenberg inequality that
\begin{align}\label{7.50}
\frac{d}{dt}\left(\|\nabla\rho\|_{L^2}+\|\nabla P\|_{L^2}\right)
& \leq  C(1+\|\nabla\mathbf{u}\|_{L^{\infty}})
\left(\|\nabla\rho\|_{L^2}+\|\nabla P\|_{L^2}+1\right)
+C\|\nabla\mathbf{b}\|_{L^\infty} \notag \\
& \leq  C(1+\|\nabla\mathbf{u}\|_{L^{\infty}})
\left(\|\nabla\rho\|_{L^2}+\|\nabla P\|_{L^2}+1\right)
+C\|\nabla\mathbf{b}\|_{L^2}^{\frac{q-2}{2q-2}}
\|\nabla^2\mathbf{b}\|_{L^q}^{\frac{q}{2q-2}} \notag \\
& \leq  C(1+\|\nabla\mathbf{u}\|_{L^{\infty}})
\left(\|\nabla\rho\|_{L^2}+\|\nabla P\|_{L^2}+1\right)
+C\|\nabla^2\mathbf{b}\|_{L^q}^{\frac{q}{2q-2}}.
\end{align}
From \eqref{2.5}, \eqref{3.1}, \eqref{6.1}, \eqref{eq7.06}, and \eqref{5.1}, we have for any $q>2$,
\begin{align}\label{6.19}
\sup_{0\leq t\leq T}\|\nabla\mathbf{u}\|_{L^q}\leq C.
\end{align}
Then, applying the standard $L^q$-estimate to \eqref{1.10}$_4$ yields
\begin{align}\label{7.6}
\|\nabla^2 \mathbf{b}\|_{L^q}&\leq C\left(\|\mathbf{b}_t\|_{L^q}
+\||\mathbf{u}||\nabla\mathbf{b}|\|_{L^q}
+\||\mathbf{b}||\nabla\mathbf{u}|\|_{L^q}\right)
\nonumber\\
&\leq C\left(\|\mathbf{b}_t\|_{H^1}
+\|\mathbf{u}\bar{x}^{-\frac{a}{2}}\|_{L^{2q}}^{\frac{q-2}{q-1}}
\|\nabla\mathbf{b}\bar{x}^{\frac{a}{2}}\|_{L^2}^{\frac{1}{q-1}}
+\|\mathbf{b}\|_{L^\infty}\|\nabla\mathbf{u}\|_{L^q}\right)\nonumber\\
&\leq C\|\nabla\mathbf{b}_t\|_{L^2}+C,
\end{align}
which combined with \eqref{7.50}, Gronwall's inequality, \eqref{7.11}, \eqref{7.13}, and the fact $\frac{q}{2q-2}\in(0,1)$ leads to
\begin{equation}\label{7.7}
\sup_{0\leq t\leq T}\left(\|\nabla \rho\|_{L^2}+\|\nabla P\|_{L^2}\right)
\leq C.
\end{equation}
Moreover, one derives from \eqref{7.3}, \eqref{3.1}, \eqref{6.1}, \eqref{7.7}, and \eqref{eq7.06} that
\begin{equation}\label{7.07}
\sup_{0\leq t\leq T}\|\nabla^2\mathbf{u}\|_{L^2}\leq C.
\end{equation}
Consequently, the desired \eqref{7.1} follows from \eqref{7.13}, \eqref{7.7}, \eqref{7.07}, and \eqref{5.1}. The proof of Lemma \ref{lem37} is finished.  \hfill $\Box$

The following higher order spatial weighted estimate on the density can be proved similarly as in \cite[Lemma 3.7]{Z2017}, and we omit the details.
\begin{lemma}\label{lem38}
Under the condition \eqref{3.1}, it holds that for any $T\in[0,T^*)$,
\begin{equation}\label{8.1}
\sup_{0\leq t\leq T}\|\rho\bar{x}^{a}\|_{H^{1}\cap W^{1,q}}\leq C.
\end{equation}
\end{lemma}

With Lemmas \ref{lem32}--\ref{lem38} at hand, we are now in a position to prove Theorem \ref{thm1.1}.

\textbf{Proof of Theorem \ref{thm1.1}.}
We argue by contradiction. Suppose that \eqref{B} were false, that is, \eqref{3.1} holds. Note that the general constant $C$ in Lemmas \ref{lem32}--\ref{lem38} is independent of $t<T^{*}$, that is, all the a priori estimates obtained in Lemmas \ref{lem32}--\ref{lem38} are uniformly bounded for any $t<T^{*}$. Hence, the function
\begin{equation*}
(\rho,\mathbf{u},P,\mathbf{b})(x,T^{*})
\triangleq\lim_{t\rightarrow T^{*}}(\rho,\mathbf{u},P,\mathbf{b})(x,t)
\end{equation*}
satisfy the initial condition \eqref{A} at $t=T^{*}$.

Furthermore, standard arguments yield that $\rho\dot{\mathbf{u}}\in C([0,T];L^2)$, which
implies $$ \rho\dot{\mathbf{u}}(x,T^\ast)=\lim_{t\rightarrow
T^\ast}\rho\dot{\mathbf{u}}\in L^2. $$
Hence, $$-\mu\Delta{\mathbf{u}}-(\lambda+\mu)\nabla\mbox{div}\mathbf{u}+\nabla P-\mathbf{b}\cdot\nabla\mathbf{b}+\frac12\nabla|\mathbf{b}|^2
|_{t=T^\ast}=\sqrt{\rho}(x,T^\ast)g(x)
$$ with $$g(x)\triangleq
\begin{cases}
\rho^{-1/2}(x,T^\ast)(\rho\dot{\mathbf{u}})(x,T^\ast),&
\mbox{for}~~x\in\{x|\rho(x,T^\ast)>0\},\\
0,&\mbox{for}~~x\in\{x|\rho(x,T^\ast)=0\},
\end{cases}
$$
satisfying $g\in L^2$ due to \eqref{7.1}.
Therefore, one can take $(\rho,\mathbf{u},P,\mathbf{b})(x,T^\ast)$ as
the initial data and extend the local
strong solution beyond $T^\ast$. This contradicts the assumption on
$T^{\ast}$.
Thus, we finish the proof of Theorem \ref{thm1.1}.
\hfill $\Box$


\end{document}